\documentclass[a4paper,10pt]{article}
\usepackage{amssymb,amsmath}
\usepackage{amssymb,enumerate}
\usepackage{graphicx}
\date{\textit{In memory of Professor Cornel Urs}}
\begin{document}
\title{Moving frames on generalized Finsler structures}
\author{S. V. Sabau\footnote{This research was partly supported by Grant-in-Aid
 for Scientific Research (C) (No. 22540097), Japan Society for the Promotion of Science.} , K. Shibuya, H. Shimada}

\maketitle
\begin{abstract} 
We study the relation between an $R$-Cartan structure $\alpha$ and an $(I,J,K)$-generalized Finsler structure $\omega$ on a 3-manifold $\Sigma$ showing the difficulty in finding a general transformation that maps $\alpha$ to $\omega$. In some particular cases, the mapping can be uniquely determined by geometrical conditions. Moreover, we are led in this way to a negative answer to our conjecture in \cite{SSS2010}.
\end{abstract}

\bigskip

{\bf MSClassification2000}: Primary 53C60, Secondary 53C20, 58A15\\
%
 


\section{Introduction}
 \quad Let  us start by recalling that a Finsler norm on a real smooth, $n$-dimensional manifold
$M$ is a function $F:TM\to \left[0,\infty \right)$ that is positive and
smooth on $\widetilde{TM}=TM\backslash\{0\}$, has the {\it homogeneity property}
$F(x,\lambda v)=\lambda F(x,v)$, for all $\lambda > 0$ and all 
$v\in T_xM$, having also the {\it strong convexity} property that the
Hessian matrix
\begin{equation*}
g_{ij}=\frac{1}{2}\frac{\partial^2 F^2}{\partial y^i\partial y^j}
\end{equation*}
is positive definite at any point $u=(x^i,y^i)\in \widetilde{TM}$.\\
\quad The fundamental function $F$ of a Finsler structure $(M,F)$ determines and it is determined by the (tangent) {\it indicatrix}, or the total space of the unit tangent bundle of $F$
\begin{equation*}
\Sigma_F:=\{u\in TM:F(u)=1\},
\end{equation*}
which is a smooth hypersurface of $TM$ such that at each $x\in M$ the {\it indicatrix at x}
\begin{equation*}
\Sigma_x:=\{v\in T_xM \ |\  F(x,v)=1\}=\Sigma_F\cap T_xM
\end{equation*}
is a smooth, closed, strictly convex hypersurface in
$T_xM$. 
   
   A Finsler structure $(M,F)$ can be therefore regarded as smooth
hypersurface $\Sigma\subset TM$ for which the canonical projection
$\pi:\Sigma\to M$ is a surjective submersion and having the property
that for each $x\in M$, the $\pi$-fiber $\Sigma_x=\pi^{-1}(x)$ is
strictly convex including the origin $O_x\in T_xM$. 

A generalization of this notion is the {\it generalized Finsler structure} introduced 
by R. Bryant (see \cite{Br1995}, \cite{Br2002}). We will consider in the following only the low dimensional case where an $(I,J,K)$-{\it generalized Finsler structure} is a 
coframing $\omega=(\omega^1,\omega^2,\omega^3)$ on a three dimensional manifold 
$\Sigma$ that satisfies the structure equations 
\begin{equation}\label{finsler_struct_eq}
\begin{split}
d\omega^1&=-I\omega^1\wedge\omega^3+\omega^2\wedge\omega^3\\
d\omega^2&=-\omega^1\wedge\omega^3\\
d\omega^3&=K\omega^1\wedge\omega^2-J\omega^1\wedge \omega^3.
\end{split}
\end{equation}
Here $I$, $J$, $K$ are smooth functions on $\Sigma$, called the invariants (in the sense of equivalence problem) of the generalized Finsler structure $
(\Sigma,\omega)$. In order to put in evidence the invariants, we will use the naming $(I,J,K)$-generalized Finsler structure.

The exterior differentials of these equations give the Bianchi type equations
\begin{equation}\label{GFS_Bianchi}
 J=I_{\omega 2}, \qquad K_{\omega 3}+KI+J_{\omega 2}=0,
\end{equation}
 where the  subscripts denote 
 the directional derivatives with respect to the coframing $\omega$, i.e.
$df=f_{\omega 1}\omega^1+f_{\omega 2}\omega^2+f_{\omega 3}\omega^3$,
for any smooth function $f$ on $\Sigma$.

By extension, one can study the generalized Finsler structure $(\Sigma,\omega)$ 
defined in this way ignoring even the existence of an underlying surface $M$.

Let us recall that, for an $(I,J,K)$-generalized Finsler structure 
$(\Sigma,\omega^1,\omega^2,\omega^3)$, we can always consider
\begin{enumerate}
\item  $ \mathcal P : \{\omega^1=0,\omega^3=0\}$ the {\it ``geodesic'' foliation} of $\Sigma$, i.e. the leaves are curves on $\Sigma$ tangent to $\hat e_2$;
\item  $ \mathcal Q : \{\omega^1=0,\omega^2=0\}$ the {\it ``indicatrix'' foliation} of $\Sigma$, i.e. the leaves are curves on $\Sigma$ tangent to $\hat e_3$,
\end{enumerate}
where $\{\hat e_1,\hat e_2,\hat e_3\}$ is the dual frame of $\{\omega^1,\omega^2,\omega^3\}$.

Any of the codimension two foliations above is called {\it amenable} if the leaf space $\Lambda$ of the foliation is a smooth surface 
such that the natural projection $\pi:\Sigma\to \Lambda$   is a smooth submersion.

By extension, an $(I,J,K)$-generalized Finsler structure $(\Sigma, \omega)$ whose ``indicatrix'' foliation $\mathcal Q$, ``geodesic'' foliation $\mathcal P$, and ``normal'' foliation 
$\mathcal R$ is amenable, is called {\it amenable}, {\it geodesically amenable} and {\it normally amenable}, respectively. Remark that in the case of a Riemannian surface, the foliations $\mathcal P$ and  $\mathcal R$ coincide up to the orientation of the Riemannian volume form, so the notions of {\it geodesically amenable} and {\it normally amenable} also coincide. Of course this is not true anymore in the case Finslerian case (\cite{SSS2010} for details).

The difference between a classical Finsler structure and a 
generalized one is global in nature, in the sense that every generalized Finsler surface structure is locally 
diffeomorphic to a classical Finsler surface structure (see \cite{Br1995}). 

The following fundamental result gives the necessary and sufficient conditions for a generalized Finsler structure to become a classical one.

{\bf Theorem 1.1.}  (\cite{Br1995})\quad {\it The necessary and sufficient condition for an $(I, J, K)$-generalized Finsler structure $(\Sigma,\omega)$ to be realizable as a classical Finsler structure on a smooth surface $M$ are
\begin{enumerate} 
\item the leaves of the indicatrix foliation $\{\omega^1=0,\ \omega^2=0\}$ to be compact;
\item it is amenable;
\item the canonical immersion $\iota:\Sigma\to TM$, given by 
$\iota(u)=\pi_{*,u}(\hat{e}_2)$, is one-to-one on each $\pi$-fiber $\Sigma_x$, where $M$ is the leaf space of the indicatrix foliation.
\end{enumerate}
 }

In order to study the differential geometry of a classical Finsler structure
$(M,F)$, it is customary to construct the pull-back bundle
$(\pi^*TM,\pi,\Sigma)$ with the $\pi$-fibers $\pi^{-1}(u)$ diffeomorphic
to $T_xM$, where $u=(x,v)\in \Sigma$ (see
\cite{BCS2000}). In general this is not a principal bundle. 

 By defining an orthonormal moving coframing on $\pi^*TM$ with respect to the
Riemannian metric on $\Sigma$ induced by the Finslerian metric $F$, the
moving equations on this frame lead to the so-called Chern
connection. This is an almost metric compatible, torsion free connection of the
vector bundle $(\pi^*TM,\pi,\Sigma)$.

\bigskip

Let us remark that the condition $I=0$ imposes $J=0$ from Bianchi identities and therefore the simplest generalized Finsler structure is a  
$(0,0,K)$-generalized Finsler structure with structure equations
\begin{equation}\label{Cartan_struct1}
\begin{split}
d\omega^1&=\omega^2\wedge\omega^3\\
d\omega^2&=-\omega^1\wedge\omega^3\\
d\omega^3&=K\omega^1\wedge\omega^2.
\end{split}
\end{equation}

This kind of structure is called a $K$-{\it Cartan structure} on the 3-manifold $\Sigma$ (see \cite{GG1995}).

We can conclude that the necessary and sufficient condition for an $(I,J,K)$-generalized Finsler structure to 
be a $K$-Cartan structure is $I= 0$.  

A $K$-Cartan structure that satisfies the conditions in Theorem 1.1 gives a Riemannian structure on the surface $M$.

A generalization of the geodesic foliation of a $K$-Cartan structure is also known (see \cite{Br1995}, \cite{Br2002}). Let us suppose that on the 3-manifold $\Sigma$ endowed with an amenable $K$-Cartan structure $\alpha$, we have a basic one-form $\beta$ for the fibration $\Sigma\to\Lambda$, 
where $\Lambda$ is the leaf space of the foliation $\{\alpha^1=0,\alpha^2=0\}$. 
Then, $ \mathcal P_\beta : \{\omega^1=0,\omega^3=\beta\}$ is called the {\it ``$\beta$-geodesic'' foliation} of $\Sigma$. Obviously, in the case $\beta=0$ we get the usual geodesic foliation. The definition can be easily extended to the non-amenable case as well. 

The geometrical meaning of $\beta$ geodesics can be thought as follows. Let us assume that the smooth surface $\Lambda$ is endowed with a Riemannian structure whose orthonormal frame bundle is $\Sigma$, and whose tautological one-forms $\alpha^1$, $\alpha^2$ together with the Levi-Civita connection form $\alpha^3$ give a coframe on $\Sigma$. Then, a smooth curve $\gamma:[a,b]\to\Lambda$ is called a $\beta$-geodesic if and only if $k_g=\beta(\dot \gamma)$, where $k_g$ is the geodesic curvature of $\gamma$ and $\dot \gamma$ is the tangent vector along $\gamma$.

\bigskip

Another particular class of special Finsler structures is the one of positive constant flag curvature $K=1$. This leads to the notion of generalized Finsler structure with $K=1$, or, an $(I,J,1)$-generalized Finsler structure on the 3-manifold $\Sigma$. The structure equations of a such structure are obtained immediately by simply substituting the curvature condition $K=1$ in \eqref{finsler_struct_eq} and \eqref{GFS_Bianchi}, respectively. Many examples of classical Finsler metrics of constant flag curvature are nowadays known, see for example \cite{BRS2004} for a complete classification of constant flag curvature Randers structures, or \cite{S2001} for other examples.

However, the existence of generalized and classical Finsler structures of constant positive flag curvature $K=1$ on surfaces was proved for the first time by R. Bryant in \cite{Br1995} using a path geometry approach. 
\bigskip

We will point out another interesting class of Finsler structures. In order to do this, let us first recall that, in the case of a classical Finsler structure 
$(M,F)$,  
the canonical parallel transport $\Phi_t:T_xM\setminus 0 \to
T_{\sigma(t)}M\setminus 0 $, 
defined by the Chern connection along a
curve $\sigma$ on $M$, is a diffeomorphism that preserves the Finslerian
length of vectors. Unlike the parallel transport on a Riemannian
manifold,  $\Phi_t$ is not a linear isometry in general. 

This unexpected fact leads to some classes of special Finsler metrics. A Finsler metric whose parallel transport is a 
linear isometry is called a 
{\it Berwald structure}, and one whose parallel transport is only a Riemannian  isometry is called a 
{\it Landsberg structure}. 

Equivalently, a Berwald metric is a Finsler metric whose Chern connection coincides with the Levi Civita 
connection of a certain Riemannian metric on $M$, in other words it is ``Riemannian-metrizable". These are the closest 
Finslerian metric to the Riemannian ones. The connection is Riemannian, while the metric is not. However, in the two 
dimensional case, any Berwald structure is Riemannian or flat locally Minkowski, i.e. there are no geometrically 
interesting Berwald surfaces. 

Landsberg structures have the property that the Riemannian volume of the Finslerian unit ball is a constant. This 
remarkable property leads to a proof of Gauss-Bonnet theorem on surfaces \cite{BCS2000} and other interesting results.

Obviously, any Berwald structure is a Landsberg one. However, there are no examples of global Landsberg 
structures that are not Berwald. The existence of Landsberg structures that are not Berwald  is one of the main open problems in modern Finsler geometry.

These notions can be easily generalized as follows. 

A {\it generalized Landsberg structure} on $\Sigma$, or an $(I,0,K)${\it -generalized Finsler structure} is a generalized Finsler 
structure $(\Sigma,\omega)$ such that $J=0$, or equivalently, $I_{\omega 2}=0$.

 Such a generalized
structure is characterized by the structure equations
\begin{equation}\label{Lands_struct_eq}
\begin{split}
d\omega^1&=-I\omega^1\wedge\omega^3+\omega^2\wedge\omega^3\\
d\omega^2&=-\omega^1\wedge\omega^3\\
d\omega^3&=K\omega^1\wedge\omega^2,
\end{split}
\end{equation}
and Bianchi identities
\begin{equation}\label{Lands_Bianchi}
\begin{split}
dI & =I_{\omega 1}\omega^1 \qquad \qquad +I_{\omega3}\omega^3 \\
dK & =K_{\omega1}\omega^1+K_{\omega2}\omega^2-KI\omega^3,
\end{split}
\end{equation}
where $I$ and $K$ are smooth functions defined on $\Sigma$.

 A {\it generalized Berwald structure} is a generalized Finsler structure characterized 
by the structure equations (\ref{Lands_struct_eq}), and
\begin{equation*}
dI \equiv 0  \quad \mod\quad \omega^3,
\end{equation*}
or, equivalently,
\begin{equation*}
I_{\omega 1}=I_{\omega 2}=0,\qquad I_{\omega 3}\neq 0.
\end{equation*}

 By means of Cartan-K\"ahler theory we have shown in \cite{SSS2010} that there exists non-trivial generalized Landsberg structures on a 3-manifold $\Sigma$. Moreover, using a path geometry approach  we have constructed locally a generalized Landsberg structure by means of a 
Riemannian structure on a surface.

\bigskip

The problems we consider in the present paper are:
\begin{enumerate}[A.]
\item {\it Clarify the relation between Cartan structures and generalized Finsler structures;
\item For given Cartan and generalized Finsler structures $\alpha$ and $\omega$, respectively, how many types of regular matrices  $A$ exist such that $\omega=A\alpha$?
\item Give a final answer to the conjecture in \cite{SSS2010} on the existence of genuine Landsberg structures on surfaces via a certain construction.
}
\end{enumerate}

The essence of the method is to start with an $R$-Cartan structure $\alpha$ on the 3-manifold $\Sigma$ and to construct by a coframe changing an $(I,J,K)$-generalized Finsler structure $\{\omega^1,\omega^2,\omega^3\}$ by means of a matrix $A=(a_{ij})$, namely 
$\omega=A\alpha$. By taking the exterior derivative of this formula and imposing conditions such that the structure equations of the $R$-Cartan structure and $(I,J,K)$-generalized Finsler structure to be satisfied, respectively, we obtain a set of differential conditions for the functions $a_{ij}$ 
in terms of its directional derivatives with respect to the coframe $\omega$. The functions satisfying these conditions can be regarded as the integral manifolds of a linear Pfaffian system. 

Using now Cartan-K\"ahler theory we can study the existence of the integral manifolds of this linear Pfaffian system. In the most general case, of course this problem is quite difficult, but in the some particular cases, for example an $(I,J,1)$-generalized Finsler structure, or, an $(I,0,K)$-generalized Finsler structure (see \cite{Br1995} and \cite{SSS2010}) one can construct explicitly the matrix $A$.

We will show in the following the difficulty of constructing such a matrix $A$ in the general case (\S2), and the simplifications appearing when one imposes conditions on the invariants $I$, $J$, $K$ and on the indicatrix foliation(\S3, \S4). Our conclusion is that subject to some geometrical conditions there are very possibilities for the matrix $A$ in the cases $K=1$ and $J=0$, respectively. 

Moreover, a Cartan-K\"ahler analysis shows that foliation conditions imposed in \S4 determine the form of the Riemannian structure on the manifold $\Lambda$ (see \S4, \S5 for notations). Namely, we ended up  with a surface of revolution $(\Lambda,g)$ whose unit sphere bundle admits an $(I,0,K)$-generalized Finsler structure. This is the same structure with the one constructed by us in \cite{SSS2010}. The present study (\S6, \S7) shows that the manifold of geodesics of the Riemannian structure $(\Lambda,g)$ is not a smooth manifold, therefore a classical Landsberg structure on a smooth surface cannot be obtained by this construction, giving in this way a negative answer to our previous conjecture in \cite{SSS2010}. 

Our findings do not rule out the existence of classical Landsberg structures on smooth surfaces, but only suggest that a different approach should be chosen.


\section{The most general case}

\quad Let us assume that we have an $R$-Cartan structure 
$\{\alpha^1,\alpha^2,\alpha^3\}$ on a 3-manifold $\Sigma$ and an $(I,J,K)$-structure $\{\omega^1,\omega^2,\omega^3\}$ on a 3-manifold $\Sigma'$. We are concerned with the problem of the existence of the linear transformation 
\begin{equation*}
A:T^*\Sigma \to T^*\Sigma'
\end{equation*}
given by $\omega=A\alpha$. Without loosing the generality we can assume that there is a diffeomorphism $\varphi:\Sigma'\to\Sigma$ and that $A$ is the matrix of the cotangent map $\varphi^*$. Therefore, formally, we identify the manifolds $\Sigma'$ and $\Sigma$.

In components, we have
\begin{equation}\label{coframe_changing_1}
\begin{pmatrix}
\omega^1\\ \omega^2\\ \omega^3
\end{pmatrix}
=\begin{pmatrix}
a_{11} & a_{12} & a_{13}\\
a_{21} & a_{22} & a_{23}\\
a_{31} & a_{32} & a_{33}
\end{pmatrix}
\begin{pmatrix}
\alpha^1\\ \alpha^2\\ \alpha^3
\end{pmatrix},
\end{equation}
where the elements $a_{ij}$ of the matrix $A$ are smooth functions on the 3-manifold $\Sigma$. Recall that an $R$-Cartan structure has the structure equations \eqref{Cartan_struct1} with the structure function $R$.

One can easily see that 
$\omega^1\wedge\omega^2\wedge\omega^3=\det(A) \ \alpha^1\wedge\alpha^2\wedge\alpha^3$.

By taking the exterior derivative of \eqref{coframe_changing_1}, it follows
\begin{equation}\label{matrix_str_eq}
d\begin{pmatrix}
\omega^1\\ \omega^2\\ \omega^3
\end{pmatrix}
=dA\wedge
\begin{pmatrix}
\alpha^1\\ \alpha^2\\ \alpha^3
\end{pmatrix}+Ad
\begin{pmatrix}
\alpha^1\\ \alpha^2\\ \alpha^3
\end{pmatrix}.
\end{equation}

Using now the structure equations of the $(I,\ J,\ K)$-structure $\omega$ and relation \eqref{coframe_changing_1}, we can express the left hand side of 
\eqref{matrix_str_eq}, with respect to the coframe $\alpha$, as follows
\begin{equation}\label{omega struct eq in terms of alpha}
\begin{split}
d\omega^1& = (-IA_{23}+A_{13})\alpha^1\wedge \alpha^2
+(-IA_{22}+A_{12})\alpha^1\wedge \alpha^3
\ +(-IA_{21}+A_{11})\alpha^2\wedge \alpha^3\\
d\omega^2&=-A_{23}\alpha^1\wedge \alpha^2
\qquad\qquad-A_{22}\alpha^1\wedge \alpha^3
\quad\qquad\qquad-A_{21}\alpha^2 \wedge\alpha^3\\
d\omega^3&=(KA_{33}-JA_{23})\alpha^1\wedge \alpha^2
+(KA_{32}-JA_{22})\alpha^1\wedge\alpha^3
+(KA_{31}-JA_{21})\alpha^2\wedge \alpha^3,
\end{split}
\end{equation}
where $A_{ij}$ is the minor of the element $a_{ij}$ in the matrix $A$, for all $i,j=1,2,3$, namely
$A_{11}=a_{22}a_{33}-a_{23}a_{32}$, $A_{12}=a_{21}a_{33}-a_{23}a_{31}$, etc.

If we denote the directional derivatives with respect to the coframe $\alpha$ by subscripts, namely
$df=f_{1}\alpha^1+f_{ 2}\alpha^2+f_{ 3}\alpha^3$,
for any smooth function $f$ on $\Sigma$, and $d a_{ij}=a_{ij\cdot 1}\alpha^1+a_{ij\cdot 2}\alpha^2+a_{ij\cdot 3}\alpha^3$ for matrices, then from \eqref{matrix_str_eq} and \eqref{omega struct eq in terms of alpha}, we obtain
the relations
\begin{equation}\label{subst1}
\begin{split}
&\begin{cases}
& a_{1 1\cdot 2}=a_{1 2\cdot 1}+a_{1 3}R +IA_{23}-A_{13}\\
& a_{2 1\cdot 2}=a_{2 2\cdot 1}+a_{2 3}R +A_{23} \\
& a_{3 1\cdot 2}=a_{3 2\cdot 1}+a_{3 3}R -KA_{33}+JA_{23} 
\end{cases}\\
&\begin{cases}
& a_{1 3\cdot 1}= a_{1 1\cdot 3}+a_{1 2} -IA_{22}+A_{12} \\
& a_{2 3\cdot 1}= a_{2 1\cdot 3}+a_{2 2} -A_{22} \\
& a_{3 3\cdot 1} = a_{3 1\cdot 3}+a_{3 2}+KA_{32}-JA_{22}
\end{cases}\\
&\begin{cases}
& a_{1 2\cdot 3} =a_{1 3\cdot 2}+a_{1 1} +IA_{21}-A_{11}\\
& a_{2 2\cdot 3}=a_{2 3\cdot 2} +a_{2 1} +A_{21}\\
& a_{3 2\cdot 3}=a_{3 3\cdot 2} +a_{3 1} - KA_{31}+JA_{21}.
\end{cases}
\end{split}
\end{equation}

Next, we will express the system of directional PDEs \eqref{subst1}, with respect to the coframe $\alpha$, using exterior differential systems (see for example \cite{IL2003} for details). 
Let us formally construct the following one forms:
\begin{equation}\label{thetas3}
\begin{split}
&\begin{cases}
& \theta_{11}=d(a_{11})-a_{11\cdot 1}\alpha^1-(a_{1 2\cdot 1}+B_{112})\alpha^2-a_{11\cdot 3}\alpha^3\\
& \theta_{21}=d(a_{21})-a_{21\cdot 1}\alpha^1-(a_{2 2\cdot 1}+B_{212})\alpha^2-a_{21\cdot 3}\alpha^3 \\
& \theta_{31}=d(a_{31})-a_{31\cdot 1}\alpha^1-(a_{3 2\cdot 1}+B_{312} )\alpha^2-a_{31\cdot 3}\alpha^3
\end{cases}\\
&\begin{cases}
& \theta_{12}=d(a_{12})-a_{12\cdot 1}\alpha^1-a_{12\cdot 2}\alpha^2-(a_{1 3\cdot 2}+B_{123})\alpha^3\\
& \theta_{22}=d(a_{22})-a_{22\cdot 1}\alpha^1-a_{22\cdot 2}\alpha^2-(a_{2 3\cdot 2} +B_{223})\alpha^3 \\
& \theta_{32}=d(a_{32})-a_{32\cdot 1}\alpha^1-a_{32\cdot 2}\alpha^2-(a_{3 3\cdot 2} +B_{323})\alpha^3
\end{cases}\\
&\begin{cases}
& \theta_{13}=d(a_{13})-(a_{1 1\cdot 3}+B_{131})\alpha^1-a_{13\cdot 2}\alpha^2-a_{13\cdot 3}\alpha^3\\
& \theta_{23}=d(a_{23})-(a_{2 1\cdot 3}+B_{231})\alpha^1-a_{23\cdot 2}\alpha^2-a_{23\cdot 3}\alpha^3 \\
& \theta_{33}=d(a_{33})-(a_{3 1\cdot 3}+B_{331})\alpha^1-a_{33\cdot 2}\alpha^2-a_{33\cdot 3}\alpha^3,
\end{cases}
\end{split}
\end{equation}
where we use the notations
\begin{align*}
 B_{112}&:= a_{1 3}R +IA_{23}-A_{13},\ & B_{212}& :=  a_{2 3}R +A_{23},\ & B_{312}& :=a_{3 3}R -KA_{33}+JA_{23}\\
 B_{123}&:= a_{1 1} +IA_{21}-A_{11},\ & B_{223}& :=  a_{2 1} +A_{21},\ & B_{323}& :=a_{3 1} - KA_{31}+JA_{21}\\
 B_{131}&:= a_{1 2} -IA_{22}+A_{12},\ & B_{231}& :=  a_{2 2} -A_{22},\ & B_{331}& :=a_{3 2}+KA_{32}-JA_{22}
\end{align*}

\bigskip

We can now formulate our setting as follows. Consider the jet space
\begin{equation*}
(x,y,p;a_{ij};a_{ij\cdot k})_{i,j,k\in\{1,2,3\}}\in J^1(3,9)=\Sigma^3\times \mathbb R^9\times \mathbb R^{27},
\end{equation*}
where $(x,y,p)$ are some local coordinates on $\Sigma$. 
Next, we consider the codimension 9 submanifold $\widetilde \Sigma^{30}$ of $J^1(3,9)$ defined by the equations \eqref{subst1}, then we can take 
the local coordinates
\begin{equation*}
(x,y,p;a_{ij};\bar{a}_{ij\cdot k})_{i,j,k\in\{1,2,3\}}\in\widetilde{\Sigma}^{30}=\Sigma^3\times \mathbb R^9\times \mathbb R^{18},
\end{equation*}
where $\bar{a}_{ij\cdot k}$ means all $a_{ij\cdot k}$, for ${i,j,k\in\{1,2,3\}}$ less the terms in the left hand side of \eqref{subst1}.

Let us consider the Pfaffian system 
\begin{equation}
\mathcal I:=\{\theta_{ij}\}_{i,j,k\in\{1,2,3\}}\subset \mathcal J:=\{\theta_{ij}; \alpha^1,\alpha^2,\alpha^3\}_{i,j,k\in\{1,2,3\}}
\subset T^*\widetilde \Sigma
\end{equation}
with independence condition $\alpha^1\wedge \alpha^2\wedge \alpha^3\neq 0$, where the one forms 
$\theta_{ij}$ are given in \eqref{thetas3}.

The integral manifolds of the linear Pfaffian system $(\mathcal I,\mathcal J)$ are the elements of the matrix $A$.

Therefore, we are going to study the existence of integral manifolds of the linear Pfaffian system 
$(\mathcal I,\mathcal J)$ defined above.

We start by computing exterior differentials of the one forms $\theta_{ij}$:
\begin{equation*}
\begin{split}
d\theta_{11}=& -d(a_{11\cdot 1})\wedge\alpha^1-d(a_{1 2\cdot 1})\wedge\alpha^2-d(a_{11\cdot 3})\wedge\alpha^3\\
& -(a_{11\cdot 1}-B_{112\cdot 3})\alpha^2\wedge\alpha^3
-(a_{1 2\cdot 1}+B_{112}) \alpha^3\wedge\alpha^1\\
&-(a_{11\cdot 3}R+B_{112\cdot 1})\alpha^1\wedge\alpha^2.
\end{split}
\end{equation*}

It follows
\begin{equation}
d\theta_{11}=\Pi_{111}\wedge\alpha^1+\Pi_{121}\wedge\alpha^2+\Pi_{113}\wedge\alpha^3,
\end{equation}
where
\begin{equation}
\begin{split}
\Pi_{111}:=& -d(a_{11\cdot 1})-(a_{1 2\cdot 1}+B_{112}) \alpha^3+(a_{11\cdot 3}R+B_{112\cdot 1})\alpha^2\\
\Pi_{121}:=& -d(a_{1 2\cdot 1})\\
\Pi_{113}:=& -d(a_{11\cdot 3})-(a_{11\cdot 1}-B_{112\cdot 3})\alpha^2.
\end{split}
\end{equation}

In the same way, we obtain
\begin{equation}
d\theta_{21}=\Pi_{211}\wedge\alpha^1+\Pi_{221}\wedge\alpha^2+\Pi_{213}\wedge\alpha^3,
\end{equation}
where
\begin{equation}
\begin{split}
\Pi_{211}:=& -d(a_{21\cdot 1})-(a_{2 2\cdot 1}+B_{212}) \alpha^3+(a_{21\cdot 3}R+B_{212\cdot 1})\alpha^2\\
\Pi_{221}:=& -d(a_{2 2\cdot 1})\\
\Pi_{213}:=& -d(a_{21\cdot 3})-(a_{21\cdot 1}-B_{212\cdot 3})\alpha^2,
\end{split}
\end{equation}
and
\begin{equation}
d\theta_{31}=\Pi_{311}\wedge\alpha^1+\Pi_{321}\wedge\alpha^2+\Pi_{313}\wedge\alpha^3,
\end{equation}
where
\begin{equation}
\begin{split}
\Pi_{311}:=& -d(a_{31\cdot 1})-(a_{3 2\cdot 1}+B_{312}) \alpha^3+(a_{31\cdot 3}R+B_{312\cdot 1})\alpha^2\\
\Pi_{321}:=& -d(a_{3 2\cdot 1})\\
\Pi_{313}:=& -d(a_{31\cdot 3})-(a_{31\cdot 1}-B_{312\cdot 3})\alpha^2.
\end{split}
\end{equation}

Therefore, we conclude
\begin{equation}
\begin{cases}
& d\theta_{11}=\Pi_{111}\wedge\alpha^1+\Pi_{121}\wedge\alpha^2+\Pi_{113}\wedge\alpha^3,\\
& d\theta_{21}=\Pi_{211}\wedge\alpha^1+\Pi_{221}\wedge\alpha^2+\Pi_{213}\wedge\alpha^3\\
& d\theta_{31}=\Pi_{311}\wedge\alpha^1+\Pi_{321}\wedge\alpha^2+\Pi_{313}\wedge\alpha^3.
\end{cases}
\end{equation}

We move now to the second group of $\theta$'s and by similar computations we 
conclude
\begin{equation}
\begin{cases}
& d\theta_{12}=\Pi_{121}\wedge\alpha^1+\Pi_{122}\wedge\alpha^2+\Pi_{132}\wedge\alpha^3\\
& d\theta_{22}=\Pi_{221}\wedge\alpha^1+\Pi_{222}\wedge\alpha^2+\Pi_{232}\wedge\alpha^3\\
& d\theta_{32}=\Pi_{321}\wedge\alpha^1+\Pi_{322}\wedge\alpha^2+\Pi_{332}\wedge\alpha^3,
\end{cases}
\end{equation}
where $\Pi_{121}$, $\Pi_{221}$, $\Pi_{321}$ are given above, and
\begin{equation}
\begin{split}
& \Pi_{122}:=-d(a_{12\cdot 2}) +(a_{12\cdot 1}+B_{123\cdot 2})\alpha^3 -(a_{1 3\cdot 2}+B_{123})R\alpha^1   \\
& \Pi_{132}:=-d(a_{1 3\cdot 2})+(a_{12\cdot 2}-B_{123\cdot 1})\alpha^1\\
& \Pi_{222}:=-d(a_{22\cdot 2})+(a_{22\cdot 1}+B_{223\cdot 2})\alpha^3-(a_{2 3\cdot 2}+B_{223})R\alpha^1               \\
& \Pi_{232}:=  -d(a_{2 3\cdot 2})+(a_{22\cdot 2}-B_{223\cdot 1})\alpha^1\\
& \Pi_{322}:=-d(a_{32\cdot 2})+(a_{32\cdot 1}+B_{323\cdot 2})\alpha^3-(a_{3 3\cdot 2}+B_{323})R\alpha^1               \\
& \Pi_{332}:=  -d(a_{3 3\cdot 2})+(a_{32\cdot 2}-B_{323\cdot 1})\alpha^1.
\end{split}
\end{equation}

Finally, we consider the third group of $\theta$'s and obtain:
\begin{equation}
\begin{cases}
& d\theta_{13}=\Pi_{113}\alpha^1+\Pi_{132}\alpha^2+\Pi_{133}\alpha^3+T_{13}\alpha^1\wedge\alpha^2\\
& d\theta_{23}=\Pi_{213}\alpha^1+\Pi_{232}\alpha^2+\Pi_{233}\alpha^3+T_{23}\alpha^1\wedge\alpha^2\\
& d\theta_{33}=\Pi_{313}\alpha^1+\Pi_{332}\alpha^2+\Pi_{333}\alpha^3+T_{33}\alpha^1\wedge\alpha^2,
\end{cases}
\end{equation}
where $\Pi_{113}$, $\Pi_{132}$, $\Pi_{213}$, $\Pi_{232}$, $\Pi_{313}$, $\Pi_{332}$ are given above, and
\begin{equation}
\begin{split}
& \Pi_{133}:=  -d(a_{13\cdot 3})-(a_{1 1\cdot 3}+B_{131}) \alpha^2
+(a_{13\cdot 2}+B_{131\cdot 3})\alpha^1          \\
& \Pi_{233}:=  -d(a_{23\cdot 3})-(a_{2 1\cdot 3}+B_{231}) \alpha^2
+(a_{23\cdot 2}+B_{231\cdot 3})\alpha^1          \\
& \Pi_{333}:=  -d(a_{33\cdot 3})-(a_{3 1\cdot 3}+B_{331}) \alpha^2
+(a_{33\cdot 2}+B_{331\cdot 3})\alpha^1          \\
& T_{13}:=-(a_{11\cdot 1}-B_{112\cdot 3})-(a_{12\cdot 2}-B_{123\cdot 1})-(a_{13\cdot 3} R-B_{131\cdot 2})\\
& T_{23}:=-(a_{21\cdot 1}-B_{212\cdot 3})-(a_{22\cdot 2}-B_{223\cdot 1})-(a_{23\cdot 3} R-B_{231\cdot 2})\\
& T_{33}:=-(a_{31\cdot 1}-B_{312\cdot 3})-(a_{32\cdot 2}-B_{323\cdot 1})-(a_{33\cdot 3} R-B_{331\cdot 2}).
\end{split}
\end{equation}

Putting all these together, we get the $\theta$'s structure equations
\begin{equation}
\begin{split}
& \begin{cases}
& d\theta_{11}=\Pi_{111}\wedge\alpha^1+\Pi_{121}\wedge\alpha^2+\Pi_{113}\wedge\alpha^3,\\
& d\theta_{12}=\Pi_{121}\wedge\alpha^1+\Pi_{122}\wedge\alpha^2+\Pi_{132}\wedge\alpha^3\\
&d\theta_{13}=\Pi_{113}\wedge\alpha^1+\Pi_{132}\wedge\alpha^2+\Pi_{133}\wedge\alpha^3
+T_{13}\alpha^1\wedge\alpha^2\\
\end{cases}
\\
& \begin{cases}
& d\theta_{21}=\Pi_{211}\wedge\alpha^1+\Pi_{221}\wedge\alpha^2+\Pi_{213}\wedge\alpha^3\\
& d\theta_{22}=\Pi_{221}\wedge\alpha^1+\Pi_{222}\wedge\alpha^2+\Pi_{232}\wedge\alpha^3\\
& d\theta_{23}=\Pi_{213}\wedge\alpha^1+\Pi_{232}\wedge\alpha^2+\Pi_{233}\wedge\alpha^3
+T_{23}\alpha^1\wedge\alpha^2\\
\end{cases}\\
& \begin{cases}
& d\theta_{31}=\Pi_{311}\wedge\alpha^1+\Pi_{321}\wedge\alpha^2+\Pi_{313}\wedge\alpha^3\\
& d\theta_{32}=\Pi_{321}\wedge\alpha^1+\Pi_{322}\wedge\alpha^2+\Pi_{332}\wedge\alpha^3\\
& d\theta_{33}=\Pi_{313}\wedge\alpha^1+\Pi_{332}\wedge\alpha^2+\Pi_{333}\wedge\alpha^3
+T_{33}\alpha^1\wedge\alpha^2.
\end{cases}
\end{split}
\end{equation}

The torsion terms are
 \begin{equation}\label{final_torsion}
\begin{split}
 T_{13}=& -(A_{13\cdot 3}+A_{11\cdot 1}-A_{12\cdot 2}) 
+I(A_{23\cdot 3}+A_{21\cdot 1}-A_{22\cdot 2})   \\
&+(I_1A_{21}-I_2A_{22}+I_3A_{23})               \\
 T_{23}=& A_{23\cdot 3}+A_{21\cdot 1}-A_{22\cdot 2} \\
 T_{33}=&-K(A_{33\cdot 3}+A_{31\cdot 1}-A_{32\cdot 2}) +J(A_{23\cdot 3}+A_{21\cdot 1}-A_{22\cdot 2})\\
&-(K_1A_{31}-K_2A_{32}+K_3A_{33})  +(J_1A_{21}-J_2A_{22}+J_3A_{23}),
\end{split}
\end{equation}
where we have used that $R_3=0$. 

{\bf Remark 2.1.}
Indeed, since $\{\alpha^1,\alpha^2,\alpha^3\}$ is an $R$-Cartan structure, from the last structure equation $d\alpha^3=R\alpha^1\wedge\alpha^2$, we get
$d^2\alpha^3=R_3\alpha^1\wedge\alpha^2\wedge\alpha^3=0$. Since the volume form $\alpha^1\wedge\alpha^2\wedge\alpha^3$ does not vanish it follows $R_3=0$. 

We can therefore conclude:

{\bf Theorem 2.1.}

{\it For a given R-Cartan structure $(\Sigma,\alpha)$ and an $(I, J, K)$-structure $(\Sigma,\omega)$, let $A:T^*\Sigma\to T^*\Sigma$ be the linear transformation whose matrix elements $(a_{ij})$ are solutions of the PDE
\eqref{subst1}. Let $(\mathcal I,\mathcal J)=(\theta^{ij},\alpha^i)$ be the corresponding Pfaffian system on the 30 dimensional manifold $\widetilde \Sigma$.

Then,  $(\mathcal I,\mathcal J)$ has non absorbable torsion, namely $T_{13}$, $T_{23}$, $T_{33}$ given by 
\eqref{final_torsion}.
}

We remark here that 
\begin{equation}\label{the system T=0}
T_{13}=0,\quad  T_{23}=0, \quad T_{33}=0
\end{equation}
is the necessary condition for the Pfaffian system $(\mathcal I,\mathcal J)$ to have solution. Following Cartan-K\"ahler algorithm (see for example \cite{IL2003}), we should restrict the Pfaffian system to the submanifold 
$$\Xi:=\{u\in\widetilde\Sigma\ |\ T_{13}(u)=0,\  T_{23}(u)=0, \ T_{33}(u)=0\}\subset \widetilde\Sigma$$
and compute again the torsion. However, the \eqref{the system T=0} is a very complicated nonlinear directional PDE, and even after restriction to $\Xi$, new non-vanishing torsions are expected to appear, so the problem is not tractable in this general case.

Therefore, we will consider in the following some particular cases when by putting some geometrical conditions we can actually decide the existence of such a linear transformation between cotangent spaces.

\section{The case $K=1$}

\quad We consider now the similar problem of an {\it a-priory} given $R$-Cartan structure 
$\alpha$ an $(I,J,1)$-structure $\omega$ on the same 3-manifold $\Sigma$. We are again concerned with the problem of the existence of the linear transformation 
$A:T^*\Sigma \to T^*\Sigma$ such that $\omega=A\alpha$. Recall that in this case
 $(I,J,1)$-generalized Finsler structure has the structure equations
 \begin{equation}\label{K=1_struct_eq}
\begin{split}
d\omega^1&=-I\omega^1\wedge\omega^3+\omega^2\wedge\omega^3\\
d\omega^2&=-\omega^1\wedge\omega^3\\
d\omega^3&=\omega^1\wedge\omega^2-J\omega^1\wedge \omega^3
\end{split}
\end{equation}
and the Bianchi identities
\begin{equation}\label{K=1_Bianchi}
  J=I_{\omega 2}, \qquad  I+J_{\omega2}=0,
\end{equation}
where the covariant derivatives are with respect to the $(I,J,1)$-structure $(\Sigma,\omega)$.

Let us remark that the most important feature of the case $K=1$ is the invariance of the two forms $\omega^1\wedge\omega^3$ and $(\omega^1)^2+(\omega^3)^2$ 
under the flow of $\hat{e}_2$.

Indeed, a straightforward computation shows that
\begin{equation}\label{2-forms omega Lie deriv}
\begin{split}
& \mathcal L_{\hat{e}_2}[(\omega^1)^2+(\omega^3)^2]=0,\qquad
\mathcal L_{\hat{e}_2} (\omega^1\wedge\omega^3)=0,\\
& \mathcal L_{\hat{e}_2}(I\omega^1+J\omega^3)=0,\qquad \quad \mathcal L_{\hat{e}_2}(I^2+J^2)=0.
\end{split}
\end{equation}

The geometrical meaning of these relations is quite clear.

If $U\subset\Sigma$ is an open set where the geodesic foliation $\{\omega^1=0,\omega^3=0\}$ is amenable (this is always possible locally), i.e. 
$\Lambda_U:=U_{\slash \{\omega^1=0,\omega^3=0\}}$ is a differential manifold and $l_U:U\to \Lambda_U$ is a smooth submersion, then there exists
\begin{enumerate}
\item a quadratic form $d\sigma^2$ on $\Lambda_U$ such that $l_U^*(d\sigma^2)=(\omega^1)^2+(\omega^3)^2$;
\item a two form $\Omega$ on $\Lambda_U$ such that $l_U^*(\Omega)=\omega^1\wedge \omega^3$;
\item a one form $\beta$ on $\Lambda_U$ such that $l_U^*(\beta)=I\omega^1+J\omega^3$.
\end{enumerate}

The indicatrix foliation of this Riemannian structure coincides with the geodesic foliation of the $(I,J,1)$-structure. 

In particular, if the $(I,J,1)$-structure $(\Sigma,\omega)$ is geodesically amenable, i.e. 
$\Lambda:=\Sigma_{\slash \{\omega^1=0,\omega^3=0\}}$ is a differential manifold and $l:\Sigma\to \Lambda$ is a smooth submersion, then on $\Lambda$ there exists a canonical Riemannian metric $g$, with area form $\Omega$ and a one form $\beta$ whose norm with respect to the Riemannian metric $g$ is $I^2+J^2$. 

On the other hand, we have given an $R$-Cartan structure $(\Sigma,\alpha)$ on the same manifold $\Sigma$ whose indicatrix foliation is $\{\alpha^1=0,\alpha^2=0\}$. Locally, if $V$ is a small neighborhood on $\Sigma$, where this foliation is amenable, then $\Lambda_V:=V_{\slash \{\alpha^1=0,\alpha^2=0\}}$ is a differential manifold and $l_V:V\to \Lambda_V$ is a smooth submersion. Moreover, $\Lambda_V$ is naturally endowed with 
\begin{enumerate}
\item a quadratic form $d\bar\sigma^2$ , such that $l_V^*(d\bar\sigma^2)=(\alpha^1)^2+(\alpha^2)^2$;
\item a two form $\bar\Omega$, such that $l_V^*(\bar\Omega)=\alpha^1\wedge \alpha^2$.
\end{enumerate}

In the case when the Cartan structure $(\Sigma,\alpha)$ is amenable, $\bar\Lambda:=\Sigma_{\slash \{\alpha^1=0,\alpha^2=0\}}$ is a differential manifold and $\bar l:\Sigma\to \bar\Lambda$ is a smooth submersion
and on $\bar \Lambda$ there exists a canonical Riemannian metric $\bar g$, with area form $\bar \Omega$.

Due to all these it is natural to impose the following 

{\bf Foliation Condition:} 

{\it The $(I,J,1)$-structure's geodesic foliation $\{\omega^1=0,\omega^3=0\}$ coincides to the indicatrix foliation of the given $R$-Cartan structure $\{\alpha^1=0,\alpha^2=0\}$.
}

This condition has the following natural

{\bf Corollary 3.1.}

{\it The $(I,J,1)$-structure $(\Sigma,\omega)$ is geodesically amenable if and only if the $R$-Cartan structure 
$(\Sigma,\alpha)$ is amenable.}

Hereafter we consider the case of an $R$-Cartan structure $(\Sigma,\omega)$ and an $(I,J,1)$-structure
 $(\Sigma,\alpha)$ under the above foliation condition. In terms of the transformation 
 $A:T^*\Sigma \to T^*\Sigma$, this implies $a_{13}=0$, $a_{33}=0$, therefore we have
\begin{equation}\label{coframe_changing_K=1}
\begin{pmatrix}
\omega^1\\ \omega^2\\ \omega^3
\end{pmatrix}
=\begin{pmatrix}
a_{11} & a_{12} & 0 \\
a_{21} & a_{22} & a_{23}\\
a_{31} & a_{32} & 0
\end{pmatrix}
\begin{pmatrix}
\alpha^1\\ \alpha^2\\ \alpha^3
\end{pmatrix}.
\end{equation}

 Another useful geometrical remark is that, in the case when $(\Sigma,\omega)$ is geodesically amenable, or, equivalently, $(\Sigma,\alpha)$ is amenable, on the 2-manifold 
 $\Lambda=\Sigma_{\slash \{\omega^1=0,\omega^3=0\}}=\Sigma_{\slash \{\alpha^1=0,\alpha^2=0\}}$
 we have actually two different quadratic forms $d\sigma^2$ and $d\bar\sigma^2$ such that 
 $l^*(d\sigma^2)=(\omega^1)^2+(\omega^3)^2$ and
 $l^*(d\bar\sigma^2)=(\alpha^1)^2+(\alpha^2)^2$, respectively. 
 In the case of amenability, they induce two different Riemannian metrics on the 2-manifold $\Lambda$. 
 
 \bigskip
 
 A natural choice is given by the following choice.
 
 {\bf Identification Condition:}
 
 {\it The Cartan structure $\{\omega^1,\omega^3\}$ coincides with the Cartan structure $\{\alpha^1,\alpha^2\}$, namely,
 \begin{equation}
 \omega^1=\alpha^1\qquad \omega^3=\alpha^2.
 \end{equation}
 }
 
 In this case, we obtain the coframe changing
\begin{equation*}
\begin{pmatrix}
\omega^1\\ \omega^2\\ \omega^3
\end{pmatrix}
=\begin{pmatrix}
1 & 0 & 0 \\
a_{21} & a_{22} & a_{23}\\
0 & 1 & 0
\end{pmatrix}
\begin{pmatrix}
\alpha^1\\ \alpha^2\\ \alpha^3
\end{pmatrix},
\end{equation*}
 and from the linear independence condition of the one forms $\alpha$ and $\omega$, namely, 
 \begin{equation*}
 \omega^1\wedge\omega^2\wedge\omega^3=-a_{23}\ \alpha^1\wedge\alpha^2\wedge\alpha^3\neq 0
 \end{equation*}
 it follows $a_{23}\neq 0$.
 
The first structure equation of \eqref{K=1_struct_eq} implies now 
\begin{equation*}
a_{23}=-1,\qquad 
a_{21}= I.
\end{equation*}
 
 The third structure equation of \eqref{K=1_struct_eq} implies further
 \begin{equation*}
 a_{22} = J,
 \end{equation*}
 so we obtain the coframe changing
 \begin{equation}\label{coframe_changing_K=1 identification case}
\begin{pmatrix}
\omega^1\\ \omega^2\\ \omega^3
\end{pmatrix}
=\begin{pmatrix}
1 & 0 & 0 \\
I & J & -1\\
0 & 1 & 0
\end{pmatrix}
\begin{pmatrix}
\alpha^1\\ \alpha^2\\ \alpha^3
\end{pmatrix}.
\end{equation}

Similarly, the second structure equation in  \eqref{K=1_struct_eq} leads us to the condition
\begin{equation}\label{curv cond for K=1, identity case}
R=-I_{2}+J_{1}+1,
\end{equation}
where the subscripts denote directional derivative with respect to the coframe $\alpha$, or, equivalently,
\begin{equation*}
R=-(J^{2}+I_{\omega3})+J_{\omega1}-I^2+1.
\end{equation*}

 Hence, we can conclude
 
 {\bf Proposition 3.2.} {\it  Let $\alpha$ and $\omega$ be an $R$-Cartan and an $(I,J,1)$-generalized Finsler structure on the 
 3-manifold $\Sigma$, respectively. Then, there is a unique  linear transformation 
$A:T^*\Sigma \to T^*\Sigma$, such that $\omega=A\alpha$, subject to the Foliation and Identification conditions above. The matrix of A in the coframes $\omega$ and $\alpha$ is given in  \eqref{coframe_changing_K=1 identification case}. In this case, the R-Cartan structure
  must obey the condition 
 \eqref{curv cond for K=1, identity case}.
 }

 \bigskip
 
{\bf  Remark.} Since we are in two dimensions and amenability assumed, these Riemannian metrics must be conformal, and therefore their tautological forms $\{\alpha^1,\alpha^2\}$ and $\{\omega^1,\omega^3\}$ must be conformal as well, respectively. 
 Here we regard the one forms 
 $\{\omega^1,\omega^3\}$ and $\{\alpha^1,\alpha^2\}$ as the tautological one forms of the corresponding Riemannian structures on the same 3-manifold $\Sigma$, respectively, up to a local diffeomorphism. Indeed, the orthogonal frame bundles of the Riemannian structures $(\Lambda,g)$ and $(\Lambda,\bar g)$ are different, but there is always possible to define local diffeomorphisms between $\Sigma$ and the total  spaces of these orthonormal frame bundles (see \cite{SSS2010} for details).
 
Therefore, without losing the generality, we can always assume the following

{\bf Conformal Equivalence}
\begin{equation}\label{conformal equiv}
\begin{pmatrix}
\omega^1\\ \omega^3
\end{pmatrix}
=m
\begin{pmatrix}
\alpha^1\\ \alpha^2
\end{pmatrix},
\end{equation}
 where $m$ is a positive function defined on $\Sigma$. From the structure equations of the $R$-Cartan structure it follows immediately that $m_{3}=0$. 
 
 Let us point out this conformal equivalence holds in the non-amenable case as well. Indeed, it can be seen from the structure equations of the $(I,J,1)$-structure that $\{\omega^1,\omega^3\}$ is a Cartan structure as well (pay attention to the fact that the third one form of this Cartan structure is different from $\omega^2$). In this case, under the foliation condition, the induced quadratic forms live on the 2-orbifold $\Lambda$ defined similarly as in the amenable case. In this case relation \eqref{conformal equiv} simply says that the two Cartan structures 
 $\{\omega^1,\omega^3\}$ and $\{\alpha^1,\alpha^2\}$ belongs to the same conformal class of Cartan structures 
 (see \cite{GG1995} for details).
 
 {\bf Remark.} One can also easily see that for a smooth function $v:\Sigma \to \mathbb{R}$,
 if $\{\alpha^1,\alpha^2,\alpha^3\}$ is an $R$-Cartan structure $\Sigma$, and $v$ a nowhere zero function, then $\{v\alpha^1,v\alpha^2,\widetilde{\alpha}^3\}$ is an $\widetilde{R}$-Cartan structure, where 
  \begin{equation}\label{Hodge_forms}
 \widetilde{\alpha}^3:=\alpha^3-*d(\log {v}),\qquad 
 v^2\widetilde{R}=R-\Delta_\alpha(\log {v}),
\end{equation}
 provided $v_{ 3}=0$.  Here $*$ and $\Delta_\alpha$ are {\it the Hodge operator} and {\it the Laplacian} of the Cartan structure $(\alpha^1, \alpha^2)$, namely
\begin{equation}
*d(\log {v})=-\frac{v_2}{v}\alpha^1+\frac{v_1}{v}\alpha^2\qquad
 \Delta_\alpha(\log {v})= \frac{1}{v^2}\Bigl[(v_{11}+v_{22})v-(v_1^2+v_2^2)     \Bigr] . 
\end{equation}
 
 The {\it conformal class}  of an $R$-Cartan structure $\{\alpha^1,\alpha^2,\alpha^3\}$ is the collection of all coframes  
 $\{v\alpha^1,v\alpha^2,\widetilde{\alpha}^3\}$, where $v$ is a smooth positive function on $\Sigma$ provided $v_3=0$.
 
 \bigskip
 
 Therefore, we obtain 
 \begin{equation}\label{coframe_changing_K=1_2}
\begin{pmatrix}
\omega^1\\ \omega^2\\ \omega^3
\end{pmatrix}
=\begin{pmatrix}
m & 0 & 0 \\
a_{21} & a_{22} & a_{23}\\
0 & m & 0
\end{pmatrix}
\begin{pmatrix}
\alpha^1\\ \alpha^2\\ \alpha^3
\end{pmatrix},
\end{equation}
 and hence
 \begin{equation*}
 \omega^1\wedge\omega^2\wedge\omega^3=-m^2a_{23}\ \alpha^1\wedge\alpha^2\wedge\alpha^3.
 \end{equation*}

The first structure equation of the $(I,J,1)$-structure gives
\begin{equation}\label{rel 1}
\begin{split}
& m_{ 3}=0\\
& a_{23}=-1\\
& a_{ 21}= (-m_{2}+Im^2)/m.
\end{split}
\end{equation}

Remark that first relation was also obtained from geometrical considerations above as well. In this case, we obtain
\begin{equation*}
\begin{pmatrix}
\omega^1\\ \omega^2\\ \omega^3
\end{pmatrix}
=\begin{pmatrix}
m & 0 & 0 \\
-\frac{m_{2}}{m}+Im & a_{22} & -1\\
0 & m & 0
\end{pmatrix}
\begin{pmatrix}
\alpha^1\\ \alpha^2\\ \alpha^3
\end{pmatrix}.
\end{equation*}

The third structure equation of $\omega$ implies
\begin{equation}
a_{22} = Jm+\frac{m_{1}}{m},
\end{equation}
 and therefore we obtain
 \begin{equation}\label{coframe_changing_K=1_3}
\begin{pmatrix}
\omega^1\\ \omega^2\\ \omega^3
\end{pmatrix}
=\begin{pmatrix}
m & 0 & 0 \\
-\frac{m_{2}}{m}+Im &  \frac{m_{1}}{m}+Jm & -1\\
0 & m & 0
\end{pmatrix}
\begin{pmatrix}
\alpha^1\\ \alpha^2\\ \alpha^3
\end{pmatrix}.
\end{equation}

Remark that we can also write formally
\begin{equation*}
\omega^2=\varphi+*d\log m-\alpha^3,
\end{equation*}
where we put $\varphi:=m(I\alpha^1+J\alpha^2)$. The Hodge operator has the usual form $*d\log m=-\frac{m_{2}}{m}\alpha^1+
\frac{m_{1}}{m}\alpha^2$.

The second structure equation \eqref{K=1_struct_eq} leads us to the following structure equation
\begin{equation}\label{phi struct eq}
d\varphi=(R-m^2)\alpha^1\wedge\alpha^2-d(*d\log m).
\end{equation}

We can conclude 

{\bf Theorem 3.3.}

{\it  

Let $\alpha$ and $\omega$ be an $R$-Cartan and an $(I,J,1)$-generalized Finsler structure on the 
 3-manifold $\Sigma$, respectively, such that the contact one-forms $\omega^1$, $\omega^3$ belong to the conformal class of Cartan structure $\alpha$. Then, the linear transformation 
$A:T^*\Sigma \to T^*\Sigma$, such that $\omega=A\alpha$, subject to the Foliation condition, is given by \eqref{coframe_changing_K=1_3}, where m is the conformal factor.

 Moreover, the structure function $R$ has to satisfy
 \begin{equation}
 R\alpha^1\wedge\alpha^2=d\varphi+m^2\alpha^1\wedge\alpha^2+d(*d\log m),
 \end{equation}
 for the one form $\varphi=m(I\alpha^1+J\alpha^2)$. 
}

Let us remark that the indicatrix foliation of the $(I,J,1)$-structure $(\Sigma,\omega)$ coincides to the 
($\varphi+*d\log m)$-geodesic foliation of the Cartan structure $(\Sigma,\alpha)$.

We can ask a 

{\bf Projective Condition:} {\it  The indicatrix foliation of the $(I,J,1)$-structure $(\Sigma,\omega)$ coincides to the geodesic foliation of the Cartan structure $(\Sigma,\alpha)$.}

This implies that we must have
\begin{equation*}
\varphi=-*d\log m,
\end{equation*}
i.e. the invariants of the $(I,J,1)$-structure $(\Sigma,\omega)$ are uniquely determined by the conformal factor $m$, namely,
\begin{equation*}
I=\frac{m_{2}}{m^2},\qquad J=-\frac{m_{1}}{m^2}.
\end{equation*}

With this extra condition, the structure equation \eqref{phi struct eq} reads
\begin{equation*}
d\varphi=-m^2\alpha^1\wedge\alpha^2+d\alpha^3-d*d\log m,
\end{equation*}
or, equivalently
\begin{equation*}
d\alpha^3=m^2\alpha^1\wedge\alpha^2,
\end{equation*}
in other words, we must have $m^2=R$.

Therefore we have

{\bf Proposition 3.4.}
 {\it Let $\alpha$ be an $R$-Cartan
structure, $R>0$  on the 
 3-manifold $\Sigma$. Then there is only one $(I,J,1)$-generalized Finsler structure on $\Sigma$ that satisfies the Foliation and the Projective conditions above, namely the one given in  \eqref{coframe_changing_K=1_3}, where  $m=\pm\sqrt{R}$ is the conformal factor. Its invariants are uniquely determined only by the structure function R, namely
 \begin{equation}
 I=\frac{1}{2}\frac{R_{2}}{R^{3/2}},\qquad  J=-\frac{1}{2}\frac{R_{1}}{R^{3/2}}.
 \end{equation} 
 }
 
 Let us recall that a Riemannian metric $g$ on a smooth manifold $\Lambda$ is called a {\it Zoll metric} if all its 
geodesics are simple closed curves of equal length. See \cite{B1978} for basics of Zoll metrics as well as 
\cite{G1976} for the abundance of Zoll metrics on $S^2$. 

 One is led in this way to the following construction method of classical Finsler structures with $K=1$ in $S^2$ found in \cite{Br2002}.
 
 {\bf Proposition 3.5.}
 {\it Let $(S^2,g)$ be a Zoll metric of positive curvature $R$ and let $\alpha^1$, $\alpha^2$ be the tautological one-forms, $\alpha^3$ the Levi-Civita connection one-form on the orthonormal sphere bundle $\Sigma$ of this Riemannian metric.
 
 Then the $(I,J,1)$-generalized Finsler structure given in Proposition 3.4 above give rise to a classical Finsler structure on the manifold of geodesics of $g$.
 }
 
The proof is straightforward by construction.


\section{The case $J=0$}


\quad We are concerned in this section with the problem of the existence of the linear transformation 
$A:T^*\Sigma \to T^*\Sigma$ such that $\omega=A\alpha$, where $\alpha$ and $\omega$ are an $R$-Cartan structure and an $(I,0,K)$-generalized Finsler structure on the same 3-manifold $\Sigma$, respectively. Recall that structure equations and Bianchi identities for $\omega$ were already given in Introduction by \eqref{Lands_struct_eq} and \eqref{Lands_Bianchi}, respectively.

One can see that it is impossible to extract any two one forms from the coframe $\omega$ to obtain a Cartan structure as we did in the previous case.

Since all forms in the coframe $\omega$ are contact forms in this case, the best thing to do is to choose an equivalent contact form from the same contact forms class. Namely, we can choose any of the following pairs of contact forms
\begin{equation}
\begin{pmatrix}
m\omega^2 \\ \omega^3
\end{pmatrix}
\quad \textrm{or}\quad
\begin{pmatrix}
\omega^2 \\m \omega^3
\end{pmatrix},
\end{equation}
where $m$ is some smooth nowhere vanishing function on $\Sigma$. More conditions will come up later. 

Any of these work, but here we choose the first one.

One can easily see that with a supplementary condition on $m$, the quadratic form 
$(m\omega^2)^2+(\omega^3)^2$ is invariant under the flow of $\hat{e}_1$. Indeed, 
\begin{equation*}
\mathcal{L}_{\hat{e}_1}((m\omega^2)^2+(\omega^3)^2)=
2\Bigl[mm_{\omega 1}(\omega^2)^2+(m^2-K)\omega^2\omega^3\Bigr].
\end{equation*}

Therefore, we get

{\bf Lemma 4.1.}

{\it Let m be a smooth non-vanishing function on $\Sigma$ satisfying
\begin{equation}\label{cond 1 on m}
m_{\omega 1}=0\qquad m^2=K,
\end{equation}
then $\mathcal{L}_{\hat{e}_1}((m\omega^2)^2+(\omega^3)^2)=0$.
}

Moreover, we have

{\bf Lemma 4.2.}

{\it Under the conditions \eqref{cond 1 on m} the pair of contact one forms
\begin{equation*}
\eta^1:=m\omega^2,\qquad \eta^2:=\omega^3, \qquad \eta^3:=m\omega^1+m_{\omega3}\omega^2
\end{equation*}
is an $R$-Cartan structure on $\Sigma$ with the structure function $R=1-\frac{m_{\omega33}}{m}$.
}

Let us return to our initial setting and denote again the coframe changing by
\begin{equation}
\begin{pmatrix}
\omega^1 \\ m\omega^2\\\omega^3
\end{pmatrix}
=\begin{pmatrix}\label{J=0 and foliation condition}
a_{11} & a_{12} & a_{13}\\
a_{21} & a_{22} & a_{23}\\
a_{31} & a_{32} & a_{33}
\end{pmatrix}
\begin{pmatrix}
\alpha^1 \\ \alpha^2\\\alpha^3
\end{pmatrix}
,
\end{equation}
where $\alpha$ and $\omega$ are arbitrary {\it a priory} given.

Similarly to the case $K=1$ all these suggest that we can put the 
following 

{\bf Foliation condition}: {\it Finslerian indicatrix foliation 
$\{\omega^1=0,\ \omega^2=0\}$ coincides with the geodesic foliation $\{\alpha^1=0,\ \alpha^3=0\}$ 
of the $R$-Cartan structure $(\Sigma,\alpha)$ (see also \cite{SSS2010})}. 

This foliation condition implies that the matrix $A$ reads now
\begin{equation}
A=\begin{pmatrix}\label{K=1 and foliation condition}
a_{11} & a_{12} & a_{13}\\
a_{21} & a_{22} & 0\\
a_{31} & a_{32} & 0
\end{pmatrix}.
\end{equation}

We can now consider again the

{\bf Identification condition:} {\it The Cartan structure $\{m\omega^2,\omega^3\}$ coincides with the R-Cartan structure $\{\alpha^1,\alpha^2\}$, i.e. 
\begin{equation*}
\alpha^1=m\omega^2,\qquad \alpha^2=\omega^3,
\end{equation*}
provided \eqref{cond 1 on m} is satisfied.}

In this case, we have
\begin{equation}
\begin{pmatrix}
\omega^1 \\ \omega^2\\\omega^3
\end{pmatrix}
=\begin{pmatrix}
a_{11} & a_{12} & a_{13}\\
\frac{1}{m} & 0 & 0\\
0 & 1 & 0
\end{pmatrix}
\begin{pmatrix}
\alpha^1 \\ \alpha^2\\\alpha^3
\end{pmatrix}
.
\end{equation}

Writing down the second and third structure equations \eqref{Lands_struct_eq} we are immediately lead to the supplementary conditions 
\begin{equation*}
a_{12}=0,\qquad a_{13}=\frac{1}{m},
\end{equation*}
and 
\begin{equation*}
a_{11}=\Bigl(\frac{1}{m}\Bigr)_{ 2}.
\end{equation*}

We therefore obtain the coframe change
\begin{equation}
\begin{pmatrix}
\omega^1 \\ \omega^2\\\omega^3
\end{pmatrix}
=\begin{pmatrix}\label{J=0, foliation, identity condition}
\Bigl(\frac{1}{m}\Bigr)_{ 2} & 0 & \frac{1}{m}\\
\frac{1}{m} & 0 & 0\\
0 & 1 & 0
\end{pmatrix}
\begin{pmatrix}
\alpha^1 \\ \alpha^2\\\alpha^3
\end{pmatrix}
.
\end{equation}

This is the same coframe changing as in \cite{SSS2010}. The first structure equation in \eqref{Lands_struct_eq} leads to
\begin{equation}
I=-2\frac{m_{ 2}}{m}, \qquad R=1-\frac{m_{ 22}}{m}.
\end{equation}

Let us point out that the Landsberg condition $I_{\omega 2}=0$ is equivalent to $m_{ 12}=0$. 

We also remark that the indicatrix foliation $\{\omega^1=0,\omega^2=0\}$ of the $(I,0,K)$-structure $(\Sigma,\omega)$ coincides with the geodesic foliation $\{\alpha^1=0,\alpha^3=0\}$ of the $R$-Cartan structure $(\Sigma,\alpha)$ by construction. 

Therefore, we can conclude

{\bf Proposition 4.3.} {\it  Let $\alpha$ and $\omega$ be an $R$-Cartan and an $(I,0,K)$-generalized Finsler structure on the 
 3-manifold $\Sigma$, respectively. Then, there is a unique  linear transformation 
$A:T^*\Sigma \to T^*\Sigma$, such that $\omega=A\alpha$, subject to the Foliation and Identification conditions above. The matrix of A in the coframes $\omega$ and $\alpha$ is given in  \eqref{J=0, foliation, identity condition}. In this case, the R-Cartan structure
  must obey the condition 
 \begin{equation}\label{curv_cond_alpha} 
R=1-\frac{m_{ 22}}{m},
\end{equation}
where m is a non-vanishing smooth function on $\Sigma$ such that
\begin{equation}\label{PDE for m}
m_{ 3}=0,\quad m_{ 12}=0.
\end{equation}
In this case, the invariants of $\omega$ are  
\begin{equation}
I=-2\frac{m_{ 2}}{m}, \qquad K=m^2.
\end{equation}
}

We remark that under the Foliation and Identification conditions, the structure functions $I$, $K$ of $(\Sigma,\omega)$ as well as the structure function $R$ of $(\Sigma,\alpha)$ are uniquely determined by the function $m$ that has to satisfy the PDE \eqref{PDE for m}.
In other words, we may say that for given $m$ solution of \eqref{PDE for m}, the 
$\Bigl(1-\frac{m_{22}}{m}\Bigr)$-Cartan structure $\alpha$ and $\Bigl(-2\frac{m_2}{m},0,m^2\Bigr)$-structure $\omega$ structures are uniquely determined by $m$.

The involutivity of the PDE made of relations in Proposition 4.3 can now be investigated by means of Cartan-K\"ahler Theorem.

Let us point out that our formula $m_{ 12}=0$ shows that the function $m_{ 1}$ is invariant under the ``geodesic" flow of $(\Sigma,\alpha)$. In other words, the function $m_{ 1}$ descends on the leaf space $\Sigma\slash_{\{\alpha^1=0,\alpha^3=0\}}$. 

\bigskip

We remark also that locally we can relax a little the Identification Condition by asking the conformal equivalence of these two Cartan structures, namely 

{\bf Conformal equivalence.} {\it 
There exists a positive smooth function $f$ on $\Sigma$, satisfying $f_{3}=0$, such that
\begin{equation*}
\begin{pmatrix}
m\omega^2\\\omega^3
\end{pmatrix}=
f\begin{pmatrix}
\alpha^1 \\ \alpha^2
\end{pmatrix}
.
\end{equation*}
}

It follows 
\begin{equation}
\begin{pmatrix}
\omega^1 \\ \omega^2\\\omega^3
\end{pmatrix}
=\begin{pmatrix}\label{K=1 and foliation condition}
a_{11} & a_{12} & a_{13}\\
\frac{f}{m} & 0 & 0\\
0 & f & 0
\end{pmatrix}
\begin{pmatrix}
\alpha^1 \\ \alpha^2\\\alpha^3
\end{pmatrix}
.
\end{equation}

We will investigate now supplementary conditions for the remaining $a_{ij}$ such that $(\Sigma,\alpha)$ and 
$(\Sigma,\omega)$ are Cartan and $(I,0,K)$-structure, respectively.

The structure equation $d\omega^2=\omega^3\wedge \omega^1$ leads to
\begin{equation}\label{formulas 1}
\begin{split}
& fa_{11}=\Bigl(\frac{f}{m}\Bigr)_{ 2}\\
& fa_{13}=\frac{f}{m}\\
&\Bigl(\frac{f}{m}\Bigr)_{ 3}=0.
\end{split}
\end{equation}

>From here it follows that $a_{11}$ and $a_{13}$ can be written in terms of $f$ and $m$ only, namely
\begin{equation}
a_{11}=\frac{1}{f}\Bigl(\frac{f}{m}\Bigr)_{ 2}
=\frac{f_{ 2}m-fm_{ 2}}{fm^2},\qquad a_{13}=\frac{1}{m}.
\end{equation}

The structure equation $d\omega^3=K\omega^1\wedge\omega^2$ leads to
\begin{equation}
a_{12}fm=-f_{ 1},\qquad f_3=0,
\end{equation}
and from here we get $a_{12}$ in terms of $f$ and $m$. Also, by combining with \eqref{formulas 1} we obtain $m_{ 3}=0$. It can be seen that this condition is consistent with $m_{1}=0$ obtained already.

With all these we obtain the final form of the coframing change
\begin{equation}\label{J=0 coframe change A}
\begin{pmatrix}
\omega^1 \\ \omega^2\\\omega^3
\end{pmatrix}
=\begin{pmatrix}
\frac{1}{f}\Bigl(\frac{f}{m}\Bigr)_{ 2} & -\frac{f_{ 1}}{fm} & \frac{1}{m}\\
\frac{f}{m} & 0 & 0\\
0 & f & 0
\end{pmatrix}
\begin{pmatrix}
\alpha^1 \\ \alpha^2\\\alpha^3
\end{pmatrix}
.
\end{equation}

It follows that the Bianchi equation $K_{\omega3}+KI=0$ implies
\begin{equation*}
I=-\frac{K_{\omega3}}{K}=-2\frac{m_{ 2}}{fm}.
\end{equation*}

After some long but not complicated computations, we get that $I_{\omega 2}=0$ is equivalent to the following
{\bf Landsberg condition}:
\begin{equation}\label{Landsberg cond A}
m_{ 21}=\frac{f_{ 1}m_{ 2}+f_{ 2}m_{ 1}}{f}.
\end{equation}

The remaining structure equation $d\omega^1=-I\omega^1\wedge\omega^3+\omega^2\wedge\omega^3$ leads now to 
\begin{equation}
\begin{split}
& -\Bigl(\frac{1}{f}\Bigl(\frac{f}{m}\Bigr)_{ 2}\Bigr)_{ 2}-\Bigl(\frac{f_{ 1}}{fm}\Bigr)_{ 1}+
\frac{1}{m}R=-I\Bigl(\frac{f}{m}\Bigr)_{2}+\frac{f^2}{m}\\
& \frac{1}{f}\Bigl(\frac{f}{m}\Bigr)_{ 2}+\Bigl(\frac{f_{ 1}}{fm}\Bigr)_{ 3}+\Bigl(\frac{1}{m}\Bigr)_{2}=I\frac{f}{m}\\
& \Bigl(\frac{1}{m}\Bigr)_{1}+\frac{f_{ 1}}{fm}-\Bigl(\frac{1}{f}\Bigl(\frac{f}{m}\Bigr)_{ 2}\Bigr)_{3}=0.
\end{split}
\end{equation}

It can be seen that the second and third relations are identities. However first relation leads to second order directional PDE for $f$ in terms of the curvature $R$ that should be considered together with the Landsberg condition \eqref{Landsberg cond A}. Indeed, by straightforward computation we get
 the {\bf curvature condition} 
\begin{equation}\label{curvature cond A}
R-f^2=\frac{f_{11}+f_{22}}{f}-\frac{f_{ 1}^2+f_{ 2}^2}{f^2}-\frac{f_{ 1}m_{ 1}-f_{ 2}m_{ 2}}{fm}
-\frac{m_{ 22}}{m}.
\end{equation}

We can conclude

{\bf Proposition 4.4.}

{\it 
Let $(\Sigma,\alpha)$ and $(\Sigma,\omega)$ be an R-Cartan structure and an $(I, 0, K)$-structure on the 3-manifold $\Sigma$, respectively. 

Let also $m$ and $f$ be two smooth functions on $\Sigma$ such that they are invariant under the flow of $\hat{e}_1$ and the conditions \eqref{Landsberg cond A} and \eqref{curvature cond A} are verified.

Then, the matrix $A$ given in \eqref{J=0 coframe change A} gives a coframing change between the R-Cartan structure $\alpha$ and the $(I,0,K)$-structure $\omega$.
}

\bigskip

It can be shown by the means of Cartan-K\"ahler theorem that the differential system 
\eqref{Landsberg cond A}, \eqref{curvature cond A} it is involutive (see \cite{SSS2010} for similar discussions), but the computations are too complicated to be given here.

\section{Cartan-K\"ahler Theory for an $(I,0,K)$-generalized Finsler structure}

\quad Let us consider here the problem of involutivity of the PDEs \eqref{curv_cond_alpha}, \eqref{PDE for m} for an unknown function $m$ on the $3$-dimensional manifold $\Sigma$.

In order to keep things simple, our setting is a 3-manifold $\Sigma$ endowed with an $R$-Cartan structure $\{\alpha^1,\alpha^2,\alpha^3\}$. All the subscripts of the scalars are again with respect to the coframe $\alpha$. 

 Let us observe that for any function $f$ on $\Sigma$ equipped with an $R$-Cartan structure, the following Ricci identities hold:
\begin{equation}
\begin{split}
& f_1+f_{32}-f_{23}=0, \\
& f_2+f_{13}-f_{31}=0, \\
& Rf_3+f_{21}-f_{12}=0.
\end{split}
\end{equation}

In particular, the Bianchi identity $R_3=0$ (see Remark 2.1.) we have 
\begin{equation}
\begin{split}
& R_1-R_{23}=0, \\
& R_2+R_{13}=0, \\
& R_{21}-R_{12}=0.
\end{split}
\end{equation}

We start by considering the following one forms 
\begin{equation}
\begin{split}
& \theta_1=dm-m_{1}\alpha^1-m_{2}\alpha^2\\
& \theta_2=dm_1-m_{11}\alpha^1+m_2\alpha^3\\
&\theta_3=dm_2-(1-R)m\alpha^2-m_1\alpha^3,
\end{split}
\end{equation}
namely, for a given $R$-Cartan structure $\{\alpha^1,\alpha^2,\alpha^3\}$ on the 3-manifold $\Sigma$, we consider PDEs \eqref{curv_cond_alpha}, \eqref{PDE for m} and Ricci identities $m_{13}=-m_2, m_{21}=0, m_{23}=m_1$.

Next, let us consider the linear Pfaffian $(\mathcal I,\mathcal J)$ with independence condition
\begin{equation}
\begin{split}
& \mathcal I=\{\theta_1,\theta_2,\theta_3\}\\
&\mathcal J=\{\alpha^1,\alpha^2,\alpha^3,\theta_1,\theta_2,\theta_3\}\\
& \alpha^1\wedge\alpha^2\wedge\alpha^3\neq 0
\end{split}
\end{equation}
on the 7-dimensional manifold $\Omega$ with the coframing
\begin{equation}
\{\alpha^1,\alpha^2,\alpha^3,dm,dm_1,dm_2,dm_{11}\}.
\end{equation}

The integral manifolds of $(\mathcal I,\mathcal J)$ are functions $m$ with the desired proprieties for given $R$-Cartan structure.

A straightforward computation shows that $(\mathcal I,\mathcal J)$ has non-absorbable torsion given by
\begin{equation}
\begin{split}
& d\theta_1 \equiv 0\\
& d\theta_2 \equiv (-dm_{11}-m_2R\alpha_2)\wedge\alpha^1+t_2\alpha^2\wedge\alpha^3, \qquad \bmod\ \theta_1, \theta_2, \theta_3\\
& d\theta_3 \equiv t_2\alpha^1\wedge\alpha^3+t_1\alpha^1\wedge\alpha^2
\end{split}
\end{equation}
where
\begin{equation}
t_1:=mR_1-m_1,\qquad t_2:=-m_{11}+m(1-R).
\end{equation}

Further, we restrict to the 5-dimensional submanifold $\widetilde \Omega:=\{t_1=0,\ t_2=0\}\subset\Omega$ with coframe
\begin{equation}
\{\alpha^1,\alpha^2,\alpha^3,dm,dm_2\}.
\end{equation}
Then, by a straightforward computation we have
\begin{eqnarray*}
\theta_2|_{\widetilde{\Omega}} &=& d(mR_1)-m(1-R)\alpha^1+m_2\alpha^3 \\
&\equiv & m\{(R_{11}+R_1^2+R-1)\alpha^1+(R_{12}+\frac{R_1m_2}{m})\alpha^2+(-R_{2}+\frac{m_2}{m})\alpha^3 \} \\
& & \hspace{7cm} \bmod\ \theta_1|_{\widetilde{\Omega}}, \theta_2|_{\widetilde{\Omega}}, \theta_3|_{\widetilde{\Omega}}
\end{eqnarray*}

>From $\theta_2|_{\widetilde{\Omega}}=0$, we obtain the equation
\begin{equation*}
-R_{2}+\frac{m_2}{m}=0 \iff m_2=mR_2
\end{equation*}
and the following conditions on the given function $R$:
\begin{equation*}
R_{11}+R_1^2+R-1=0,\ R_{12}+R_1R_2=0.
\end{equation*}

>From now on, we assume given $R$ satisfies above two equations. We restrict again to the 4-dimensional submanifold $\widetilde{\widetilde \Omega}:=\{m_2=mR_2\}\subset\widetilde\Omega$ with coframe
\begin{equation}
\{\alpha^1,\alpha^2,\alpha^3,dm\}
\end{equation}

It follows
\begin{eqnarray*}
\theta_3|_{\widetilde{\widetilde{\Omega}}} &=& d(mR_{2})-(1-R)m\alpha^2-mR_1\alpha^3\\
&\equiv & m(R_{22}+R_2^2+R-1)\alpha^2 \qquad \bmod\ \theta_1|_{\widetilde{\widetilde{\Omega}}}, \theta_2|_{\widetilde{\widetilde{\Omega}}}, \theta_3|_{\widetilde{\widetilde{\Omega}}}
\end{eqnarray*}

Hence, in order to solve the PDEs \eqref{curv_cond_alpha}, \eqref{PDE for m}, the function $R$ must satisfy:
\begin{equation}
R_{22}+R_2^2+R-1=0.
\end{equation}

>From now on, we assume that the $R$-Cartan structure $\alpha$ satisfies
\begin{equation}\label{EDS for R}
\left\{
\begin{split}
& dR-R_{ 1}\alpha^1-R_{2}\alpha^2=0 \\
& dR_1+(R_1^2+R-1)\alpha^1+R_1R_2\alpha^2+R_2\alpha^3 =0 \\
& dR_2+R_1R_2\alpha^1+(R_2^2+R-1)\alpha^2-R_1\alpha^3 =0.
\end{split}
\right.
\iff
\left\{
\begin{split}
& R_3=0\\
& R_{11}=1-R-R_1^2\\
& R_{12}=R_{21}=-R_1R_2\\
& R_{22}=1-R-R_2^2.
\end{split}
\right.
\end{equation}

Hence, on the 4-dimensional submanifold $\widetilde{\widetilde \Omega}:=\{m_2=mR_2\}\subset\widetilde\Omega$ with the coframe
\begin{equation}
\{\alpha^1,\alpha^2,\alpha^3,dm\}
\end{equation}
we have $\theta_2|_{\widetilde{\widetilde{\Omega}}}=\theta_3|_{\widetilde{\widetilde{\Omega}}}=0$ and 
\begin{equation}
d \theta_1|_{\widetilde{\widetilde{\Omega}}}= m(d(\log m-R)).
\end{equation}
Therefore, the Pfaffian system $\{\theta_1|_{\widetilde{\widetilde{\Omega}}}=0\}$ on $\widetilde{\widetilde \Omega}:=\{m_2=mR_2\}\subset\widetilde\Omega$ is a Frobenius system.

>From the above computations, we remark that the following relation holds good
\begin{equation}\label{exp R}
dm-mdR=0 \Longleftrightarrow m=Ce^R,
\end{equation}
where $C$ is arbitrary constant, in other words, $m$ is uniquely determined by the structure function $R$ of the Cartan structure $\alpha$.

We conclude

{\bf Proposition 5.1.} {\it For a given R-Cartan structure which satisfies 
\eqref{EDS for R}, we have one parameter family of solutions of the PDEs 
\eqref{curv_cond_alpha}, \eqref{PDE for m}.}

Taking into account of \eqref{exp R} we can also formulate

{\bf Proposition 5.2.} {\it
Let us consider an $R$-Cartan structure $\alpha$ on the 3-manifold $\Sigma$ which satisfies \eqref{EDS for R}.

Then, the coframe 
\begin{equation}
\begin{split}
&\omega^1=\frac{1}{C}e^{-R}(\alpha^3-R_2\alpha^1) \\
&\omega^2=\frac{1}{C}e^{-R}\alpha^1\\
& \omega^3= \alpha^2
\end{split}
\end{equation}
gives an $(I,0,K)$ generalized Finsler structure (i.e. a generalized Landsberg structure) on the 3-manifold $\Sigma$ with the invariants
\begin{equation}
I=-2R_2,\qquad K=C^2 e^{2R}.
\end{equation}
}

{\bf Remark 5.1.} The existence of an $R$-Cartan structures which satisfy \eqref{EDS for R} is guaranteed by the existence of Generalized Landsber structures (see \cite{SSS2010}) or one can
prove the existence of the $R$-Cartan structures by using the same method
(Cartan-K\"ahler theorem for the frame bundle) as in \cite{SSS2010}.


\section{The geometry of the leaf space $\Lambda$}

\quad We will study in the present section the geometry of the quotient space $\Lambda$ that is the leafs space of the codimension two foliation 
$\{\alpha^1=0,\alpha^2=0\}$ on the 3-manifold $\Sigma$. This section is based on a private communication with R. Bryant (\cite{Br2011}). 

Let us recall from the previous section that we are given an $R$-Cartan structure $\{\alpha^1,\alpha^2,\alpha^3\}$ on $\Sigma$ whose structure function $R$ satisfies
\eqref{EDS for R}.

{\bf Theorem 6.1.} 

{\it The $R$-Cartan structure $(\Sigma,\alpha)$ subject to the conditions \eqref{EDS for R} induces on the leaf space
$\Lambda=\Sigma_{\slash\{\alpha^1=0,\alpha^2=0\}}$ a family of smooth rotationally symmetric Riemannian metrics
\begin{equation}\label{g 1}
g=\eta^2(r)dr^2 + \varphi^2(r)d\theta^2,
\end{equation}
where
\begin{equation}
\eta(r)=\frac{re^{-\frac{1}{4}r^2}}{2\Phi(r)},\qquad \varphi(r)=\frac{\Phi(r)}{R_0-1},
\end{equation}
and
\begin{equation}
\Phi(r)=\Bigl[(R_0-\frac{3}{2})-(R_0-\frac{3}{2}-\frac{1}{4}r^2)e^{-\frac{1}{2}r^2}\Bigr]^\frac{1}{2}.
\end{equation}
Here $(r,\theta)\in I\times S^1$ are the induced local coordinates on $\Lambda$, and $R_0>1$ is a real constant.
}

{\bf Proof.}
It can be seen that from equations \eqref{EDS for R}, one obtains
\begin{equation}
d\Bigl(\bigl(R_1^2+R_2^2+R-\frac{3}{2}\bigr)e^{2(R-1)}\Bigr)=0,
\end{equation}
therefore, it must exists a constant $C$ such that
\begin{equation}\label{eq for R 1}
\bigl(R_1^2+R_2^2+R-\frac{3}{2}\bigr)e^{2(R-1)}=C.
\end{equation}

Observe that the last formula can be rewritten as 
\begin{equation}
(R_1^2+R_2^2)e^{2(R-1)}+f(R)=C,
\end{equation}
where 
\begin{equation}
f(R)=\bigl(R-\frac{3}{2}\bigr)e^{2(R-1)}
\end{equation}
is a smooth function of one variable. Since $(R_1^2+R_2^2)e^{2(R-1)}\geq 0$ it follows that we must have $f(R)\leq C$. A simple computation shows that $f(1)=-\frac{1}{2}$ is the minimum of $f$.
%
%

Therefore, one gets $C\geq -\frac{1}{2}$ with equality in the case $R=1$, which is not good for us and it will be therefore eliminated. 

For any fixed $C> -\frac{1}{2}$, if we denote by $(f(R_0),R_0)$ the intersection point of the graph of $f$ with the straight line $f=C$, it means that always there exists a unique constant $R_0>1$ such that $C=f(R_0)=\bigl(R_0-\frac{3}{2}\bigr)e^{2(R_0-1)}$, therefore \eqref{eq for R 1} reads now
\begin{equation}\label{eq for R 2}
\bigl(R_1^2+R_2^2+R-\frac{3}{2}\bigr)e^{2(R-1)}=\bigl(R_0-\frac{3}{2}\bigr)e^{2(R_0-1)},
\end{equation} 
for any constant $R_0>1$. Obviously we must have $R\leq R_0$.

On the other hand, one can see that
\begin{equation}
d\Bigl[ \frac{e^{(R_0-R)}}{R_1^2+R_2^2}(R_1\alpha^1+R_2\alpha^2)
\Bigr]=d\Bigl[ \frac{e^{(R_0-R)}}{R_1^2+R_2^2}(R_2\alpha^1-R_1\alpha^2)
\Bigr]=0
\end{equation}

Therefore, if we denote
\begin{equation}
dz^1=\frac{e^{(R_0-R)}}{R_1^2+R_2^2}dR,\qquad dz^2=\frac{e^{(R_0-R)}}{R_1^2+R_2^2}dR
\end{equation}
it follows
\begin{equation}
(dz^1)^2+(dz^2)^2=\frac{e^{2(R_0-R)}}{R_1^2+R_2^2}\bigl[(\alpha^1)^2+(\alpha^2)^2\bigr].
\end{equation}

Therefore, at least locally, one obtains a metric on $\Lambda$
\begin{equation}
g=(\alpha^1)^2+(\alpha^2)^2=\frac{R_1^2+R_2^2}{e^{2(R_0-R)}}\bigl[(dz^1)^2+(dz^2)^2\bigr],
\end{equation} 
by using now \eqref{eq for R 2} we get the equivalent form
\begin{equation}
\begin{split}
g & = \Bigl((R_0-\frac{3}{2})-(R-\frac{3}{2})e^{2(R-R_0)}\Bigr)\bigl[(dz^1)^2+(dz^2)^2\bigr]\\
   & = \frac{e^{2(R-R_0)}}{(R_0-\frac{3}{2})-(R-\frac{3}{2})e^{2(R-R_0)}}dR^2+
   \Bigl((R_0-\frac{3}{2})-(R-\frac{3}{2})e^{2(R-R_0)}\Bigr){(dz^2)}^2.
\end{split}
\end{equation}

One could now regard $(R,z^2)$ as the local coordinates on $\Lambda$, but it is better to use the coordinate change $(r,\theta)\mapsto(R,z^2)$, where
\begin{equation}
R=R_0-\frac{r}{4}^2,\qquad z^2=\frac{\theta}{R_0-1}.
\end{equation}

This puts the metric $g$ in the ``polar coordinates" defined on a domain $I\times S^1$
\begin{equation}
g=\eta^2(r)dr^2 + \varphi^2(r)d\theta^2,
\end{equation}
where $I$ is an interval in $\mathbb R$ and
\begin{equation}
\begin{split}
& \eta(r)=\frac{1}{2}\frac{re^{-\frac{1}{4}r^2}}{\Bigr[(R_0-\frac{3}{2})-(R_0-\frac{3}{2}
-\frac{1}{4}r^2)e^{-\frac{1}{2}r^2}\Bigr]^\frac{1}{2}} \\
& \varphi(r) =\frac{\Bigl[(R_0-\frac{3}{2})-(R_0-\frac{3}{2}-\frac{1}{4}r^2)e^{-\frac{1}{2}r^2}\Bigr]^\frac{1}{2}}{R_0-1},
\end{split}
\end{equation}
or, equivalently,
\begin{equation}
\eta(r)=\frac{re^{-\frac{1}{4}r^2}}{2\Phi(r)},\qquad \varphi(r)=\frac{\Phi(r)}{R_0-1},
\end{equation}
where we put
\begin{equation}
\Phi(r)=\Bigl[(R_0-\frac{3}{2})-(R_0-\frac{3}{2}-\frac{1}{4}r^2)e^{-\frac{1}{2}r^2}\Bigr]^\frac{1}{2}.
\end{equation}

One can easily see that this is a {\it smooth rotational metric} even at the origin $r=0$ and that the curvature $R$ reaches its maximum there.

\begin{flushright}
$\Box$
\end{flushright}

{\bf Remark 6.1.} 
\begin{enumerate}
\item Conversely, for any $R_0>1$, formula \eqref{g 1} defines a metric on a region of the plane where $(R_0-\frac{3}{2})-(R_0-\frac{3}{2}-\frac{1}{4}r^2)e^{-\frac{1}{2}r^2}>0$ that satisfies the desired conditions.

\item The value of the parameter $R_0>1$ influences the geometry of the Riemannian metric $g$, indeed
\begin{enumerate}
\item for $R_0\geq\frac{3}{2}$, $(r,\theta)\in \mathbb R\times S^1$, i.e.  $(r,\theta)$ are global coordinates, in other words $(\Lambda,g)$ is a smooth Riemannian manifold;
	\begin{enumerate}
	\item when $R_0 > \frac{3}{2}$, 
	the metric \eqref{g 1} is incomplete (because the $r$-curves have finite length), and can not be extended because $$\lim_{t\to\infty}k=-\infty;$$
	\item when $R_0=\frac{3}{2}$, the metric gets the form
	\begin{equation}\label{g for R_0=3/2}
	g=dr^2+r^2e^{-\frac{r^2}{2}}d\theta^2
	\end{equation}
	which is a complete Riemannian metric;
	\end{enumerate}
\item for $1<R_0<\frac{3}{2}$ the condition 
	$(R_0-\frac{3}{2})-(R_0-\frac{3}{2}-\frac{1}{4}r^2)e^{-\frac{1}{2}r^2}\geq 0$ holds for $|r|\leq T(R_0)$ for a certain function $T(R_0)>0$. 
	
	This metric is not complete at $t=\pm T(R_0)$, but the length of the parallels goes to 0 as $t\to\pm T(R_0)$. Thus, the metric ``completes" to a singular 2-sphere, with one singularity at the ``minimum curvature pole". This will be an actual singularity if the cone angle is not $2\pi$, and this will depend on the value of $R_0$.
\end{enumerate}
\end{enumerate}

Let us point out that the prime integral $m_1$ of the geodesic flow of $(\Lambda,g)$ constructed in Section 4 reads now
$m_1=R_1e^R$ and therefore it gives a first coordinate on the geodesic space $M=\Sigma/_{\{\alpha^1=0,\alpha^3=0\}}$. Since the space $M$ is the leaf space of a codimension two foliation, it should have another coordinate, but to find it explicitly on $\Sigma$, we need a second prime integral of the geodesic flow of $(\Lambda,g)$, and this is a difficult problem to be studied in the future. Instead of doing this, we will induce local coordinates on $M$ from the tangent bundle of $\Lambda$ in the next section.

\section{The geometry of the surface of revolution $(\Lambda,g)$}

Let us consider the complete surface of revolution with the metric \eqref{g for R_0=3/2}. 
We recall that the geodesic flow of a surface of revolution is completely integrable. 

In polar coordinates, we have
\begin{equation}
g_{11}=1,\quad g_{12}=g_{21}=0,\quad g_{22}=h^2(r)=r^2e^{-\frac{r^2}{2}}.
\end{equation}
The sectional curvature of this metric is
\begin{equation}\label{sect_curvat}
R=\frac{3}{2}-\frac{1}{4}r^2\leq \frac{3}{2},
\end{equation}
and non-vanishing Christoffel symbols
\begin{equation}
\gamma^1_{22}=-hh'=\frac{r(r^2-2)}{2}e^{-\frac{1}{2}r^2},\quad 
\gamma^2_{12}=\frac{h'}{h}=\gamma^2_{21}=-\frac{1}{2}\frac{r^2-2}{r}.
\end{equation}

It follows that a curve $\gamma:[a,b] \to \Lambda$, with $\gamma(t)=(r(t),\theta(t))$ is a geodesic on $\Lambda$ if and only if it satisfies the following geodesic equations:
\begin{equation}
\begin{split}
&r''+\frac{r(r^2-2)}{2}e^{-\frac{1}{2}r^2}(\theta')^2=0\\
&\theta''-\frac{r^2-2}{r}r'\theta'=0.
\end{split}
\end{equation}


A straightforward computation shows that 
\begin{equation}
\frac{d}{dt}\Bigl[  h^2(r)\theta'
\Bigr]_{\gamma(t)}=0,
\end{equation}
i.e. the function
\begin{equation}
F=h^2(r)\theta'=r^2e^{-\frac{r^2}{2}}\theta'
\end{equation}
is constant on the geodesic lines, and thus provides a prime integral of the geodesic flow of $(\Lambda,g)$. This is called the {\it Clairaut integral} (see for example \cite{SST2003}). 

A natural question to ask is if our prime integral $m_1$ is different from the usual Clairaut integral. With the notations above, we have

{\bf Proposition 7.1.} {\it The prime integral $R_1e^R$ coincides with the Clairaut integral F up to multiplication with a negative constant.}

{\it Proof.} Let us denote by $f_1$, $f_2$ the $g$-orthonormal frame on $\Lambda$, then the relation with the natural frame 
$\dfrac{\partial}{\partial r}$, $\dfrac{\partial}{\partial \theta}$ can be written as
\begin{equation}
\begin{split}
&  f_1=   \frac{1}{\sqrt {|g|}}\Bigl( \frac{\partial\sqrt E}{\partial \zeta^2}  \frac{\partial}{\partial r}- \frac{\partial\sqrt E}{\partial \zeta^1}  \frac{\partial}{\partial \theta}\Bigr)
=:m^1 \frac{\partial}{\partial r}+m^2 \frac{\partial}{\partial \theta} \\
& f_2=\frac{1}{\sqrt E}\Bigl(\zeta^1\frac{\partial}{\partial r}+\zeta^2\frac{\partial}{\partial \theta}\Bigr)=:l^1 \frac{\partial}{\partial r}+l^2 \frac{\partial}{\partial \theta},
\end{split}
\end{equation}
where $E=(\zeta^1)^2+h^2(r)(\zeta^2)^2$ is the energy function of the Riemannian structure $(\Lambda,g)$ and $|g|=h^2$ is the determinant of the matrix $(g_{ij})$.
Here we use $(r,\theta,\zeta^1,\zeta^2)\in T\Lambda$ to denote the natural coordinates on the tangent space of $(\Lambda,g)$.

Their horizontal lifts to the unit sphere bundle $\Sigma$ of $(\Lambda,g)$ are
\begin{equation}
  \hat f_1=   m^1 \frac{\delta}{\delta r}+m^2 \frac{\delta}{\delta \theta},\qquad
 \hat f_2=l^1 \frac{\delta}{\delta r}+l^2 \frac{\delta}{\delta \theta},
\end{equation}
where $\dfrac{\delta}{\delta r}$, $\dfrac{\delta}{\delta \theta}$ is the adapted basis of the horizontal subspace of $T\Lambda$ with respect to the Levi-Civita connection, i.e.
\begin{equation}
\begin{split}
& \frac{\delta}{\delta r}=\frac{\partial}{\partial r}-\gamma^2_{12}\zeta^2\frac{\partial}{\partial \zeta^2}=\frac{\partial}{\partial r}-\frac{h'}{h}\zeta^2\frac
{\partial}{\partial \zeta^2}\\
& \frac{\delta}{\delta \theta}=\frac{\partial}{\partial \theta}-\gamma^1_{22}\zeta^2\frac{\partial}{\partial \zeta^1}
-\gamma^2_{21}\zeta^1\frac{\partial}{\partial \zeta^2}
=\frac{\partial}{\partial \theta}+hh'\zeta^2\frac{\partial}{\partial \zeta^1}
-\frac{h'}{h}\zeta^1\frac{\partial}{\partial \zeta^2}.
\end{split}
\end{equation}

With these notations, we have
\begin{equation}
R_1=m^1 \frac{\delta R}{\delta r}+m^2 \frac{\delta R}{\delta \theta}=m^1\frac{\partial R}{\partial r}.
\end{equation}

We will compute now $m^1$ as follows
\begin{equation}
m^1=\frac{1}{\sqrt {|g|}}\frac{\partial\sqrt E}{\partial \zeta^2}=\frac{1}{2E\sqrt {|g|}}\frac{\partial E}{\partial \zeta^2}=h(r)\frac{\zeta^2}{E},
\end{equation}
and taking into account of \eqref{sect_curvat} it follows that on $\Sigma$ we have
\begin{equation}
R_1=-rh(r)\zeta^2.
\end{equation}

On the other hand, let us recall that the canonical lift of a geodesic $\gamma$ of $(\Lambda,g)$ to the unit sphere bundle $\Sigma=\mathcal F(\Lambda,g)$ gives the geodesic flow of $(\Lambda,g)$ on $\Sigma$. Its tangent vector is
\begin{equation}
\hat f_2=\hat T(t)=\dot\gamma^i(t)\frac{\delta}{\delta x^i}=r'\frac{\delta}{\delta r}+\theta'\frac{\delta}{\delta \theta},
\end{equation} 
and therefore, along the lines of the geodesic flow we have $\theta'=d\theta(\hat f_2)=l^2\equiv \zeta^2$ on $\Sigma$.

Finally, observing that
\begin{equation}
e^R=e^{\frac{3}{2}-\frac{1}{4}r^2}=e^\frac{3}{2}\frac{h}{r}
\end{equation}
it follows
\begin{equation}
R_1e^R=-e^\frac{3}{2}h^2\zeta^2,
\end{equation}
i.e., on the geodesic flow of $(\Lambda,g)$ we have $R_1e^R=-e^\frac{3}{2}F$, and the Proposition is proved.
\begin{flushright}
$\Box$
\end{flushright}

Finally, let us remark that we can regard the space $M$ as a submanifold of the tangent space $T\Lambda$. On $T\Lambda$ we have already two functionally independent prime integrals of the geodesic flow of $(\Lambda,g)$, namely, the Clairaut integral $F=-e^{-\frac{3}{2}}R_1e^R$, regarded now as function on $T\Lambda$, and the energy function $E$ of the Riemannian structure $(\Lambda,g)$. 

Therefore, we obtain in this way induced local coordinates $(F,E)$ on $M$ from $T\Lambda$, i.e.,
$\pi_M:T\Lambda\to M$, $(r,\theta,\zeta^1,\zeta^2)\mapsto (F,E)$. Pay attention to the fact that $M$ is not a smooth manifold because $(\Lambda,g)$ cannot have all geodesics closed and of same length.

\bigskip

{\bf Acknowledgments.}  We would like to express our special thanks to R. Bryant for many useful discussions on the subject, as well as for his invaluable help in understanding the geometry of the space $\Lambda$ that led to Section 6.



\begin{thebibliography}{[Br et al 1991]}

\bibitem[BCS2000]{BCS2000}
     {Bao},~D., {Chern},~S.S., {Shen},~Z.,
     {\it An Introduction to Riemann Finsler Geometry}, Springer, GTM 200, 
2000.


\bibitem[BRS2004]{BRS2004}
	{Bao},~D., Robles,~C.,  {Shen}, ~Z.,
	{\it Zermelo navigation on Riemannian manifolds},
J. Diff. Geom., {\bf Vol. 66} (2004), 391--449.

\bibitem[B1978]{B1978}
	{Besse},~A.,
	{\it Manifolds all of whose Geodesics are Closed}, Springer-Verlag, 1978.

\bibitem[Br1996]{Br1995}
	       {Bryant},~R.,
{\it Finsler structures on the 2-sphere satisfying $K=1$}, Finsler Geometry,
Contemporary Mathematics {\bf 196} (1996), 27--41.

\bibitem[Br2002]{Br2002}
 {Bryant},~R.,
{\it Some remarks on Finsler manifolds with constant
flag curvature, Houston Journal of Mathematics}, {\bf 28}, no.2,
(2002), 221--262.

\bibitem[Br2011]{Br2011}
{Bryant},~R., {\it Private communication}, January 2011.

\bibitem[GG1995]{GG1995}
		{Geiges},~H., {Gonzalo},~J.,
{\it Contact Geometry and Complex Surfaces},
Invent. Math., {\bf 121}(1995), 147--209.


\bibitem[GG1997]{GG1997}
		{Geiges}, ~H., {Gonzalo},~J.,
{\it Contact circles on 3-manifolds},
J. Diff. Geom., {\bf 46} (1997), 236--286.

\bibitem[GG2002]{GG2002}
		{Geiges}, ~H., {Gonzalo},~J.,
{\it Moduli of contact circles},
J. Reine Angew. Math., {\bf 551}(2002), 41--85. 


\bibitem[G1976]{G1976}
		{Guillemin},V.,
		{\it The Radon transform on Zoll surfaces},
		Advances in Math. {\bf 22} (1976), 85--119.



\bibitem[IL2003]{IL2003}
{Ivey},~Th.~A., {Landsberg},~J.~M.,
{\it Cartan for Beginners; Differential Geometry via Moving Frames and
		 Exterior Differential systems}, AMS, GSM 61, 2003.

\bibitem[SSS2010]{SSS2010}
{Sabau},~S. V., {Shibuya},~K., {Shimada},~H.,
{\it On the existence of generalized unicorns on surfaces}, Diff. Geom. and its Appl., {\bf 28} (2010), 406--435.

\bibitem[S2001]{S2001}
	{Shen}, ~Z.,
	{\it Lectures on Finsler Geometry}, World Scientific, 2001. 


\bibitem[SST2003]{SST2003}
	{Shiohama}, ~K., {Shioya}, ~T., {Tanaka}, ~M.,
{\it The geometry of total curvature on complete open surfaces},
Cambridge University Press, 2003.


\end{thebibliography}
\end{document}